\documentclass{lms}
\usepackage{rotating}
\newtheorem{theorem}{Theorem}[section]

\newtheorem{keylemma}[theorem]{Key Lemma}

\renewcommand{\tt}{\mathfrak{t}}

\newcommand{\cue}{\mathbb{Q}}
\newcommand{\cee}{\mathbb{C}}

\newcommand{\W}{{\mathcal W}}
\renewcommand{\S}{{\mathcal S}}

\renewcommand{\L}{{\mathcal L}}

\newcommand{\N}{{\mathcal N}}

\newcommand{\cptwobar}{\overline{\cee P}^2}
\newcommand{\Mod}{\mbox{Mod}}

\usepackage{graphicx,amsfonts}

\title{Monodromy Substitutions and Rational Blowdowns}

\author{Hisaaki Endo, Thomas E. Mark and Jeremy Van Horn-Morris}

\classno{57R65, 20F38, 57R17}

\extraline{The first author was partially supported by Grand-in-Aid for Scientific Research (C), 
No. 21540079, Japan Society for the Promotion of Science. The second author was partially supported by NSF grant DMS-0905380.}

\begin{document}
\maketitle

\begin{abstract}We introduce several new families of relations in the mapping class groups of planar surfaces, each equating two products of right-handed Dehn twists. The interest of these relations lies in their geometric interpretation in terms of rational blowdowns of 4-manifolds, specifically via monodromy substitution in Lefschetz fibrations. The simplest example is the lantern relation, already shown by the first author and Gurtas \cite{EG} to correspond to rational blowdown along a $-4$ sphere; here we give relations that extend that result to realize the ``generalized'' rational blowdowns of Fintushel-Stern \cite{FS} and Park \cite{park} by monodromy subsitution, as well as several of the families of rational blowdowns discovered by Stipsicz-Szab\'o-Wahl \cite{SSW}.
\end{abstract}

\section{Introduction}

It has been known for some time that a rational blowdown of a smooth 4-manifold can be performed symplectically. Moreover, symplectic 4-manifolds are well-known to correspond to Lefschetz pencils or, after suitable blowup, to Lefschetz fibrations. Our aim here is to combine these two ideas to show that under certain circumstances one can perform a rational blowdown on a 4-manifold equipped with a Lefschetz fibration, and preserve the fibration structure. 

Recall that a Lefschetz fibration on a closed smooth 4-manifold $X$ is a smooth map $\pi: X\to S^2$ that is a fiber bundle projection away from finitely many singular points. Near the singular points $\pi$ is required to appear in appropriate oriented local complex coordinates as $\pi(z_1,z_2) = z_1z_2$. The Lefschetz fibration is said to have genus $g$ if the typical fiber of $\pi$ is a smooth surface $\Sigma$ of genus $g$. A theorem of Gompf (\cite{GS}, Theorem 10.2.18) shows that under mild hypotheses the total space of a Lefschetz fibration admits a symplectic structure. 

The monodromy of $\pi$ around a critical value is well-known to be isotopic to a (right-handed) Dehn twist around a simple closed curve in $\Sigma$, the {\it vanishing cycle} associated to the critical point. Arranging the critical values $p_1,\ldots, p_n\in S^2$ in a cyclic order, the Lefschetz fibration is determined up to isomorphism by the {\it global monodromy} given by the sequence of corresponding Dehn twists. Observe that the composition of all these Dehn twists is isotopic to the identity (it is the  sequence itself that determines $X$, not the composition). Hoping it does not lead to confusion, we often refer to this global monodromy as an (unreduced) {\it word} in right-handed Dehn twists, though the monodromy is determined only up to cyclic permutation of the twists and simultaneous conjugation by some diffeomorphism of the fiber $\Sigma$. We write $X = X(w)$ for the Lefschetz fibration determined by the word $w$.

For a surface $F$ let $\Mod(F)$ (or $\Mod(F,\partial F)$) denote the mapping class group of $F$, the group of isotopy classes of orientation-preserving diffeomorphisms of $F$ (or of $F$, fixing the boundary pointwise). Suppose that $w_1$ and $w_2$ are two words in right-handed Dehn twists with the property that $w_1 = w_2$ in $\Mod(\Sigma)$, and suppose $X = X(w_1\cdot w')$ is a Lefschetz fibration with fiber $\Sigma$ and global monodromy given by $w_1\cdot w'$ for some word $w'$ That is to say, $w_1$ and $w_2$ are distinct factorizations of a given element of $\Mod(\Sigma)$ into products of Dehn twists, and one of these factorizations appears in the monodromy word associated to $X$. By a {\it monodromy substitution}, we mean the operation 
\[
X = X(w_1\cdot w') \mapsto X' = X(w_2\cdot w')
\]
of replacing $X$ by the Lefschetz fibration $X'$ having global monodromy $w_2 \cdot w'$. Observe that the products $w_1w'$ and $w_2w'$ are equal in $\Mod(\Sigma)$, so that the boundaries of $X$ and $X'$ are naturally diffeomorphic. In particular, if $w_1w'$ is isotopic to the identity so that $X$ is a closed Lefschetz fibration, the same is true of $X'$. Generally $X$ and $X'$ will be different smooth 4-manifolds: indeed, if $w_1$ and $w_2$ involve a different number of Dehn twists, then $X$ and $X'$ have differing Euler characteristics. However, $X'$ does still possess the structure of a Lefschetz fibration (with the same fiber) and hence is symplectic.

\bigskip\noindent
{\bf Examples:} 1) If $\gamma$ is a curve that bounds a disk on $\Sigma$, then the corresponding Dehn twist is isotopic to the identity. Writing $\gamma$ also for the right-handed Dehn twist around $\gamma$, suppose $X = X(\gamma\cdot w')$ for some word $w'$. The result of monodromy substitution based on the relation $\gamma \sim 1$ in $\Mod(\Sigma)$ gives rise to the fibration $X(w')$ with monodromy $w'$. As a smooth $4$-manifold, $X(w')$ is obtained from $X(\gamma\cdot w')$ by contracting a sphere of self-intersection $-1$, i.e., an ordinary blowdown operation.

2) The {\it lantern relation} states that $abcd = xyz$, where the letters indicate (twists around) the curves shown in Figure \ref{lanternfig} below---thinking of the 3-holed disk in that figure as a subsurface of $\Sigma$. It was shown in \cite{EG} that in the corresponding monodromy substitution $X = X(abcd\cdot w')\mapsto X' = X(xyz\cdot w')$, the resulting manifold $X'$ is obtained from $X$ by cutting out a neighborhood of a sphere with normal bundle of degree $-4$ and replacing it by a rational homology ball: a {\it rational blowdown}.  

\bigskip
It is natural, particularly in light of the second example, to ask whether there are other relations in appropriate mapping class groups that correspond to other rational blowdowns under monodromy substitution. By a rational blowdown, we will mean the operation of replacing the neighborhood of a configuration of spheres in a smooth 4-manifold, intersecting according to some connected plumbing graph, by a rational ball having the same oriented boundary. There are many examples of plumbed 4-manifolds whose boundary also bounds a rational ball, the first having been studied in this context by Fintushel and Stern \cite{FS}. A generalization of their examples to a family of linear plumbing graphs $\L$ was considered by Park \cite{park}, and several more families were discovered by Stipsicz-Szab\'o-Wahl \cite{SSW}. Our main result shows that many of these rational blowdowns can be realized by monodromy substitution. 

\begin{theorem}\label{mainthm} Let $G$ be any plumbing graph among the families $\mathcal L$, $\mathcal W$, or $\mathcal N$ (see below). Then there exists a planar surface $S$ and a relation of the form $w_1 = w_2$ in $\Mod (S,\partial S)$, where $w_1$ and $w_2$ are words in right-handed Dehn twists, with the following property. Suppose that $X$ is a Lefschetz fibration with fiber $\Sigma$ contaning $S$ as a subsurface, and global monodromy of the form $w_1 \cdot w'$. Then $X = X(w_1\cdot w')$ contains an embedded copy of the plumbing of disk bundles over spheres determined by $G$, and furthermore the Lefschetz fibration $X' = X(w_2\cdot w')$ given by monodromy substitution is diffeomorphic to the rational blowdown of $X$ along this plumbing.
\end{theorem}

The families $\mathcal W$ and $\mathcal N$ consist of graphs $\Gamma_{p,q,r}$ and $\Delta_{p,q,r}$, respectively, for arbitrary nonnegative integers $p,q,r$; they are shown in Figure \ref{plumbingsfig}.
The family $\mathcal L$ contains linear plumbing graphs $C_{p,q}$ indexed by a pair of relatively prime integers $p > q> 0$, with vertices weighted by the integers $-b_1,\ldots, -b_k$. The $b_i$ are the continued fraction coefficients of $p^2/(pq-1)$:
\[
\frac{p}{q} = b_1 - \frac{1}{b_2 - \frac{1}{\cdots -\frac{1}{b_k}}}.
\]

\begin{figure}[t]
\includegraphics{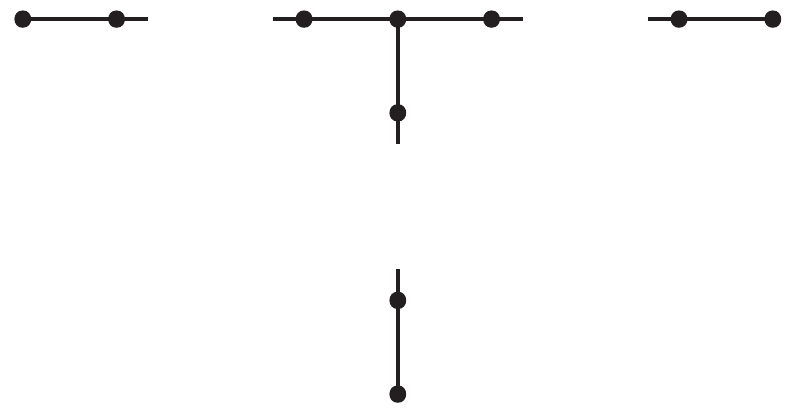}
\put(-250,100){$-(p+3)$}
\put(-200,108){$\underbrace{\hspace{.9in}}_{q}$}
\put(-174,111){$\cdots$}
\put(-120,120){$-4$}
\put(-93,108){$\underbrace{\hspace{.9in}}_r$}
\put(-65,111){$\cdots$}
\put(-20,100){$-(q+3)$}
\put(-118,57){$\left.\begin{array}{c}\vspace*{.68in} \\ \end{array}\right\}$}
\put(-115,55){$\vdots$}
\put(-95,57){\scriptsize $p$}
\put(-105,5){$-(r+3)$}
\put(-220,-15){(a): plumbing graph $\Gamma_{p,q,r}$ in the $\mathcal W$ family}
\vspace{2ex}\\
\includegraphics{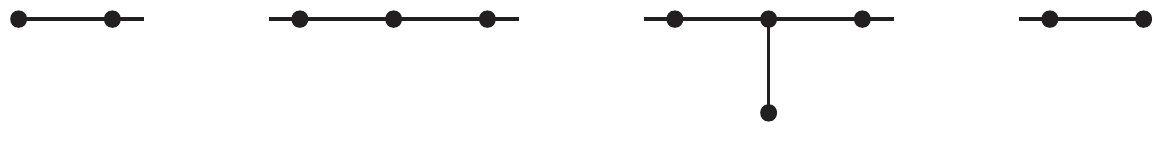}
\put(-355,42){$-(r+3)$}
\put(-308,30){$\underbrace{\hspace{.9in}}_{p-1}$}
\put(-200,30){$\underbrace{\hspace{.9in}}_q$}
\put(-93,30){$\underbrace{\hspace{.9in}}_{r}$}
\put(-283,33){$\cdots$}
\put(-175,33){$\cdots$}
\put(-66,33){$\cdots$}
\put(-230,42){$-3$}
\put(-120,42){$-3$}
\put(-25,42){$-(q+4)$}
\put(-130,-5){$-(p+2)$}
\put(-270,-20){(b): graph $\Delta_{p,q,r}$ ($p\geq 1$) in the $\mathcal N$ family}
\vspace{2ex}
\\
\includegraphics{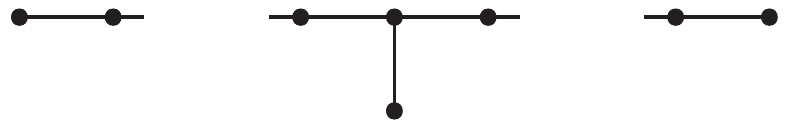}
\put(-250,42){$-(r+4)$}
\put(-125,42){$-3$}
\put(-25,42){$-(q+4)$}
\put(-200,30){$\underbrace{\hspace{.9in}}_q$}
\put(-92,30){$\underbrace{\hspace{.9in}}_r$}
\put(-175,33){$\cdots$}
\put(-66,33){$\cdots$}
\put(-150,-20){(c): graph $\Delta_{0,q,r}$}
\caption{\label{plumbingsfig} Rational blowdown plumbing graphs. Unlabeled vertices carry weight $-2$.}
\end{figure}
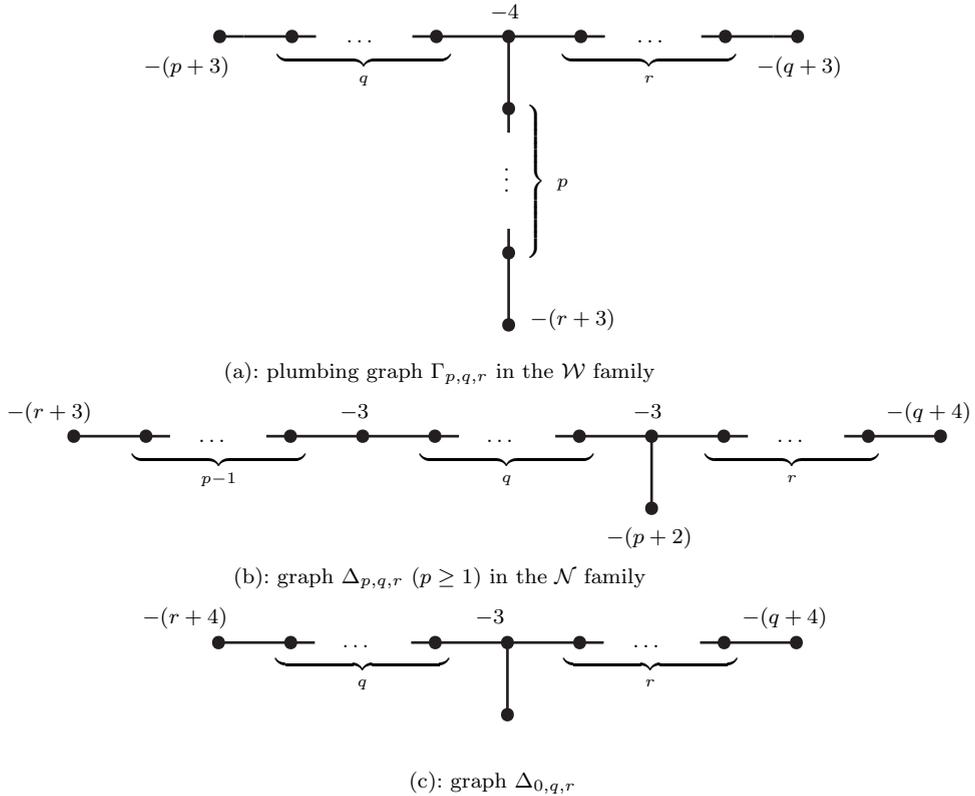
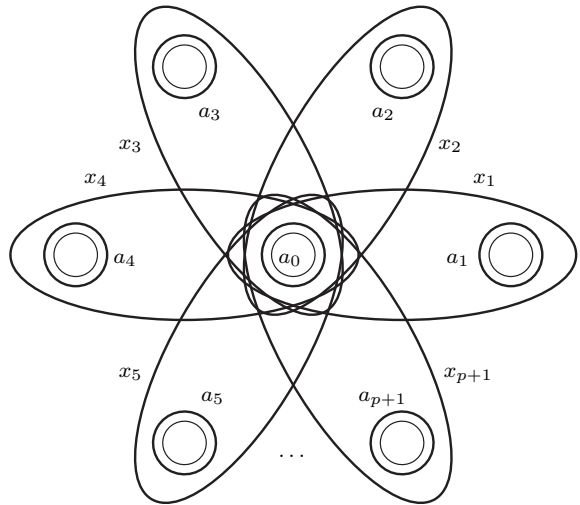
\begin{figure}[h]
\includegraphics[width=3in]{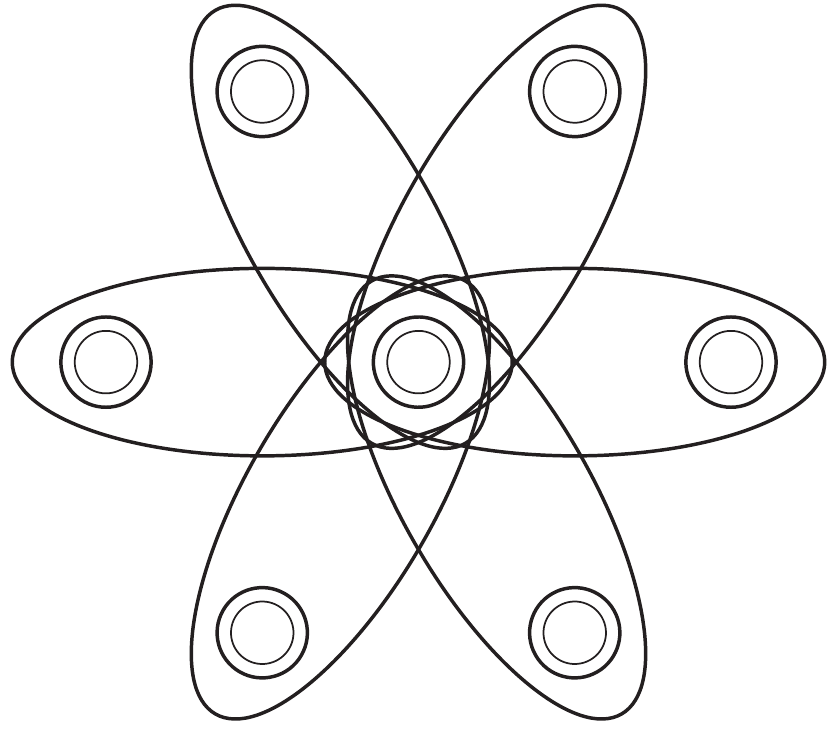}
\put(-50,93){$a_1$}
\put(-113,93){$a_0$}
\put(-78,148){$a_2$}
\put(-143,148){$a_3$}
\put(-175,93){$a_4$}
\put(-142,41){$a_5$}
\put(-83,41){$a_{p+1}$}
\put(-113,20){$\cdots$}
\put(-40,123){$x_1$}
\put(-53,136){$x_2$}
\put(-173,136){$x_3$}
\put(-186,123){$x_4$}
\put(-173,50){$x_5$}
\put(-51,50){$x_{p+1}$}
\caption{\label{daisyfig} The daisy relation. The figure is drawn on a sphere with ${p+2}$ holes (indicated by the lighter circles). The relation states that the composition of right-handed Dehn twists around the curves $x_1,\ldots, x_{p+1}$ (in that order) is isotopic to the product of twists around the boundary-parallel curves $a_0,\ldots ,a_{p+1}$, with the twist around $a_0$ raised to the power $p-1$. A monodromy substitution given by replacing the product of $a$'s by the product of $x$'s results in a rational blowdown along a linear plumbing graph $C_{p,1}$.}
\end{figure}

A simple example of a relation arising in the theorem is  the one corresponding to Fintushel-Stern's original construction: a linear plumbing graph $G_p$ ($p\geq 2$, corresponding to $C_{p,1}$) whose vertices are weighted by the (negatives of the) coefficients of the continued fraction expansion of $p^2/(p-1)$ (that is, by $\{-(p+2), -2, -2, \ldots, -2\}$ where $-2$ occurs $p-2$ times). The corresponding relation is in the mapping class group of a sphere with $p + 2$ holes; the case $p=2$ is just the lantern relation. The relation for general $p$ is shown in Figure \ref{daisyfig}, and is called a ``daisy relation'' for obvious reasons. We remark that this relation has appeared also in the work of Plamenevskaya and the third author \cite{OJ}.

We close this introduction with a couple of remarks. First, Theorem \ref{mainthm} is entirely local: it proceeds by finding a Lefschetz fibration $N(G)\to D^2$ on the neighborhood $N(G)$ of a plumbing in the listed families, having a planar surface as fiber. Then we apply a relation to write the monodromy of the corresponding open book on $\partial N(G)$ as a different product of Dehn twists, and the corresponding different Lefschetz fibration is a rational ball with diffeomorphic boundary. In particular, there is no need to assume that $X$ is a closed Lefschetz fibration to begin with.  On the other hand, the hypotheses of Theorem \ref{mainthm} appear to be somewhat restrictive, in the sense that if one encounters a (symplectic) manifold containing a plumbing among the listed families, it is not obvious whether a global Lefschetz fibration structure exists having a monodromy of the required form. However it can be shown that the local Lefschetz fibration on $N(G)$ can be extended to $X$ as a {\it broken} Lefschetz fibration (c.f. \cite{AK}, \cite{baykur}, \cite{GK}, \cite{lekili}). It is possible, moreover, that the purely local structures introduced here can be used to show that rational blowdowns (or more general operations) obtained by monodromy substitution may be performed symplectically. We leave this as a future project.

Finally, there are several families of rational blowdowns (notably the family $\mathcal M$ of \cite{SSW}) not covered by Theorem \ref{mainthm}. It is an interesting question to decide whether such blowdowns can also be recovered by substitution techniques. On the other hand, the methods of this paper yield many more relations in planar mapping class groups (or in mapping class groups of surfaces of higher genus). So far, only those relations in Theorem \ref{mainthm} have been given an ``interesting'' geometric interpretation, but it is natural to hope for new examples of operations on 4-manifolds arising from these ideas. 

{\bf Organization:} The relations needed for our monodromy substitutions are derived in the following section. On some level, the proof of Theorem \ref{mainthm} is entirely elementary once the relations are in hand, and requires only some basic checking. However, in section \ref{diagramsec} we produce Kirby diagrams for the rational balls appearing in our rational blowdowns, and verify that they are diffeomorphic to those found in \cite{FS}, \cite{park}, and \cite{SSW}. In the final section we construct a family of closed Lefschetz fibrations over $S^2$ to which a monodromy substitution corresponding to a daisy relation may be applied, yielding minimal symplectic manifolds homeomorphic but not diffeomorphic to $\#(4n^2 - 4n +1)\cee P^2 \#(12n^2 - 2n)\overline{\cee P}^2$ for each $n\geq 2$.

{\bf Acknowledgements:} The first author is grateful to Kouichi Yasui for helpful discussions; we also thank Andr\'as Stipsicz for his interest and encouragement. 

\section{Monodromy relations}

\begin{figure}[h]
\includegraphics[width=3in]{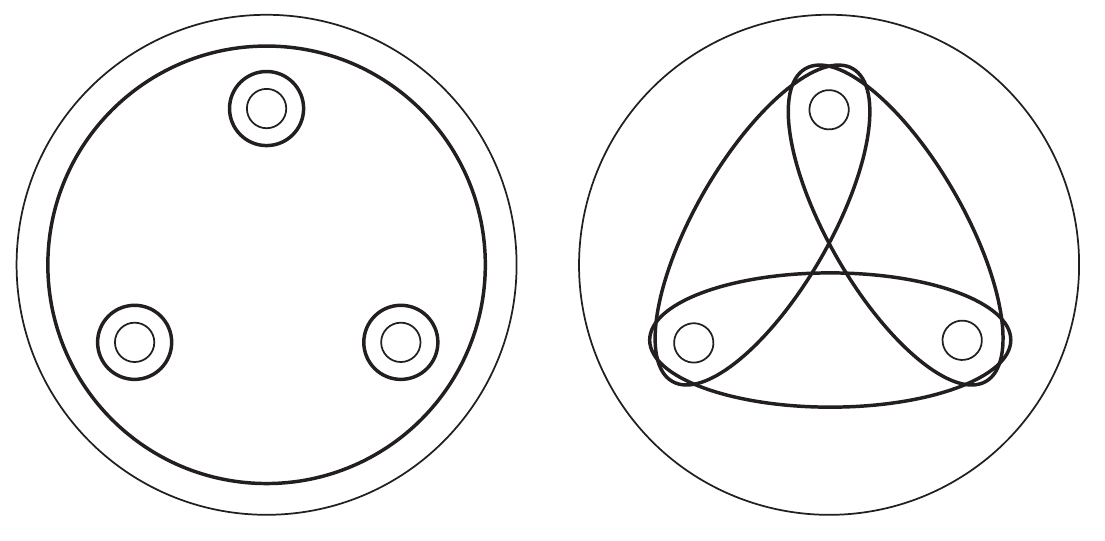}
\put(-170,70){$a$}
\put(-150,50){$b$}
\put(-185,50){$c$}
\put(-200,70){$d$}
\put(-88,70){$x$}
\put(-25,70){$y$}
\put(-55,17){$z$}
\caption{\label{lanternfig}Two copies of the 4-holed sphere $S_4$. Lighter circles are boundary components; darker circles are curves to twist around.}
\end{figure}
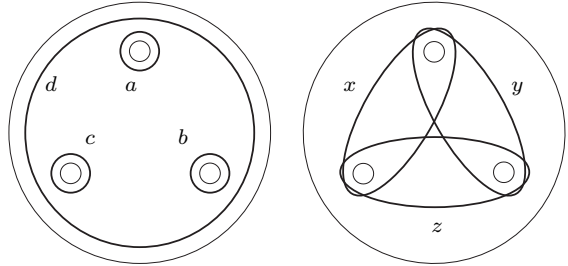

All of our relations may be derived from the basic lantern relation (indeed, the planar mapping class group may be given a presentation in which the lantern is the only relation, aside from commutation relations: see \cite{margalit}). Recall that if $S_4$ is a sphere with four holes and $a,b,c,d,x,y,z$ are the curves indicated in Figure \ref{lanternfig}, then the lantern relation states that 
\[
abcd= xyz.
\]
Here we are sloppy and use the same symbol for (the isotopy class of) a curve and (the isotopy class of) the right-handed Dehn twist about that curve, and in the above relation we use group multiplication: that is, in the product $xy$, the twist around $x$ precedes the twist around $y$. Note that the surface $S_4$ of Figure \ref{lanternfig} may also be a subsurface of a more general surface.

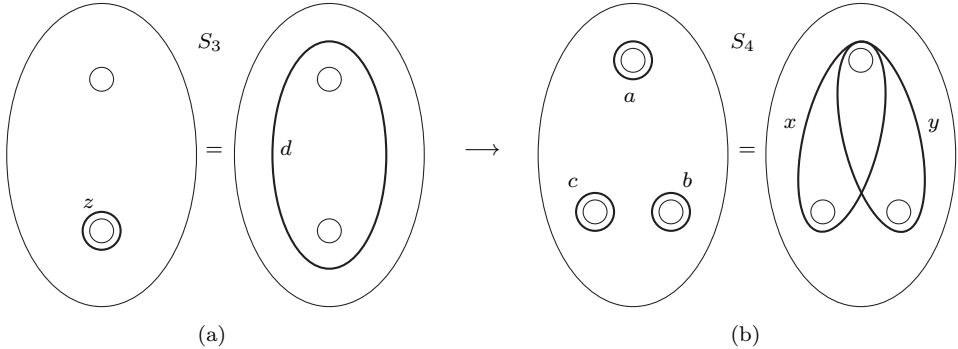
\begin{figure}
\includegraphics[width=5in]{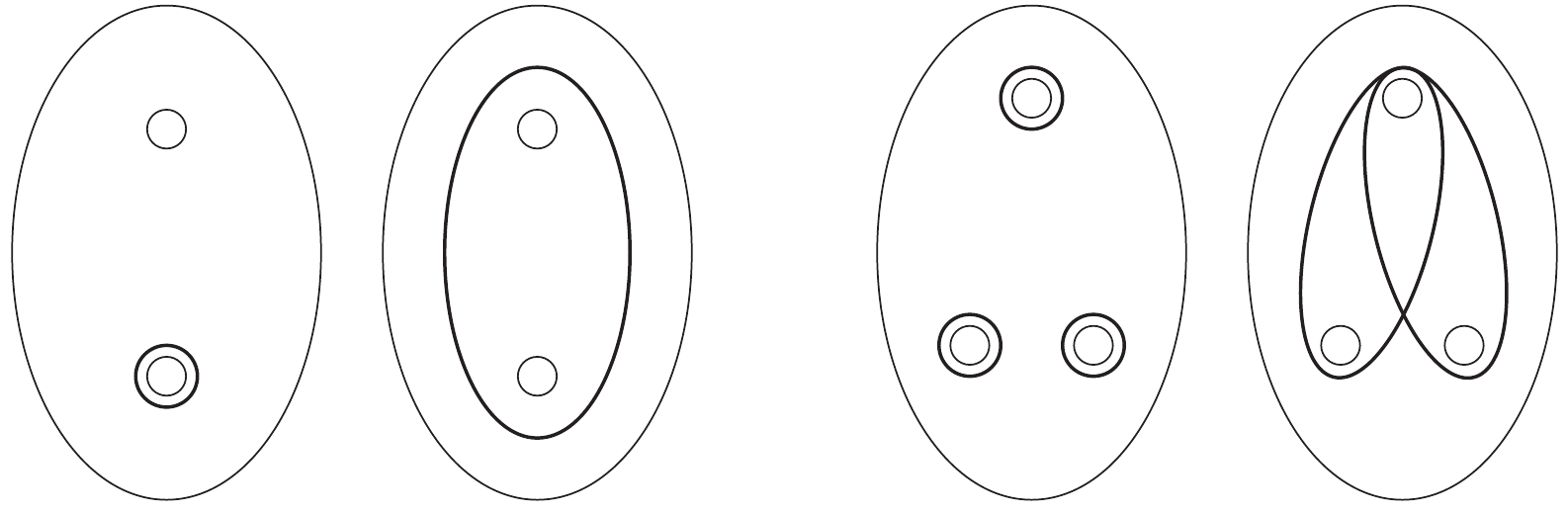}
\put(-330,40){$z$}
\put(-256,60){$d$}
\put(-148,48){$c$}
\put(-105,48){$b$}
\put(-127,80){$a$}
\put(-67,70){$x$}
\put(-13,70){$y$}
\put(-287,-10){(a)}
\put(-87,-10){(b)}
\put(-284,60){$=$}
\put(-84,60){$=$}
\put(-187,60){$\longrightarrow$}
\put(-287,100){$S_3$}
\put(-87,100){$S_4$}
\caption{In (a) we see two copies of a subsurface $S_3 \subset F$, and in (b) two copies of a subsurface $S_4$ obtained by removing a disk from $S_3$. The ``$=$'' are meant to suggest relations between words in the mapping class group, where in (a) (a twist around) $z$ appears on one side of the relation and $d$ on the other. We obtain a new relation on a new surface by the operation suggested in (b).}\label{keylemmafig}
\end{figure}

Our generalizations of this relation are all based on the following simple construction. Consider a planar surface $F$  containing as a subsurface the pair-of-pants $S_3$, and let $z$ and $d$ be the boundary-parallel curves marked as in Figure \ref{keylemmafig}(a). 
We construct a new surface $F'$ with one additional boundary component by removing a small disk from $F$ (and $S_3$), the closure of which we take to lie in the interior of the region between the curve $z$ and the hole of $S_3$ it encloses. If the component of $\partial S_3$ corresponding to $z$ coincides 
with a component of $\partial F$, then we can think of $F'$ as obtained from $F$ by gluing a disk with two holes into the hole enclosed by $z$. (Intuitively, we think of this as ``splitting'' the hole enclosed by $z$ in two.) Extending by the identity across this new 2-holed disk induces a homomorphism $\varphi: \Mod(F,\partial F)\to \Mod(F',\partial F')$ that we use implicitly in the following lemma. Note that in this situation, the twist around $z$ commutes with the image of $\varphi$.
 
\begin{keylemma} \label{keylemma} Given the setup of the previous paragraph (in particular, $z$ encloses only one hole of $F$),
suppose that in the planar mapping class group $\Mod(F,\partial F)$ the relation 
\[
w_1zw_2 = w_1'dw_2'
\]
holds, for some $w_1,w_2,w_1',w_2'\in\Mod(F,\partial F)$. Assume that $a$ commutes with either $w_1$ and $w_1'$, or with $w_2$ and $w_2'$. Then in $\Mod(F',\partial F')$ we have the relation
\[
w_1abcw_2 = w_1'xyw_2'
\]
where $a,b,c,x,y$ are the curves in $F'$ marked in Figure \ref{keylemmafig}(b).
\end{keylemma}

\begin{proof} Since $z$ commutes with $w_1$ and $w_1'$ (either on $F$ or $F'$) and with $d$, the given relation shows 
\[
w_1w_2 = w_1'dz^{-1}w_2'
\]
in $\Mod(F',\partial F')$. Multiply both sides by $abc$ on the left (if $a$ commutes with $w_1$, $w_1'$) or right (in the other case), and observe that both $b$ and $c$ are central in $\Mod(F',\partial F')$. The lantern relation completes the proof.
\end{proof}

Observe that the commutativity assumption on $a$ is trivially satisfied if the hole of $S_3$ corresponding to $a$ coincides with a hole of $F$: in this case $a$ is boundary parallel and therefore central.

Our plan in this section is to inductively construct several families of relations on planar surfaces $F$ using the Key Lemma. For reasons that will become clear subsequently, we are looking for relations of the form $w = w'$ satisfying the following properties:
\begin{enumerate}
\item Both $w$ and $w'$ are words in right-handed Dehn twists.
\item $w$ consists of twists around a collection of pairwise disjoint curves.
\item If $F$ is an $n$-holed disk, then $w'$ is a product of $n$ twists around $n$ curves that span the rational first homology of $F$.
\end{enumerate}
The prototypical example of such a relation is the lantern relation itself. Observe that if relation $R_2$ is obtained from relation $R_1$ by successive applications of the Key Lemma \ref{keylemma}, and if $R_1$ satisfies (1), (2), and (3) above, then $R_2$ will satisfy (1) and (3). By careless choices of subsurface $S_3$, however, it is easy to disrupt property (2).

\subsection{Relations corresponding to family $\mathcal W$}\label{Wrelsubsec}

\begin{theorem}[($\mathcal W$-family of relations)]\label{Wrelthm} For three integers $p,q,r\geq 0$, let $F$ be a disk with $p+q+r+6$ holes, with corresponding boundary-parallel curves $a_1,\ldots, a_{p+2}$, $b_1,\ldots, b_{q+2}$, $c_1,\ldots, c_{r+2}$ (arranged on a circle and numbered clockwise for convenience). Let $d_a$, $d_b$, $d_c$ be disjoint simple closed curves respectively enclosing only the $a_i$, $b_i$ or $c_i$ (see Figure \ref{Wrelfig}), and let $d$ be a circle parallel to the outer boundary of $F$. For $i = 1,\ldots, p+2$, let $A_i$ be a ``convex'' curve enclosing exactly the holes corresponding to $c_1,\ldots, c_{r+2}$ and $a_i$. Similarly let $B_1,\ldots, B_{q+2}$ enclose all $a$'s and a single $b$, and $C_1,\ldots, C_{r+2}$ enclose all $b$'s and a single $c$. Then in $\Mod(F,\partial F)$ we have
\[
a_1\cdots a_{p+2}b_1\cdots b_{q+2}c_1\cdots c_{r+2}d_a^{q+1}d_b^{r+1}d_c^{p+1}d = A_{p+2}\cdots A_1B_{q+2}\cdots B_1C_{r+2}\cdots C_1.
\]
\end{theorem}
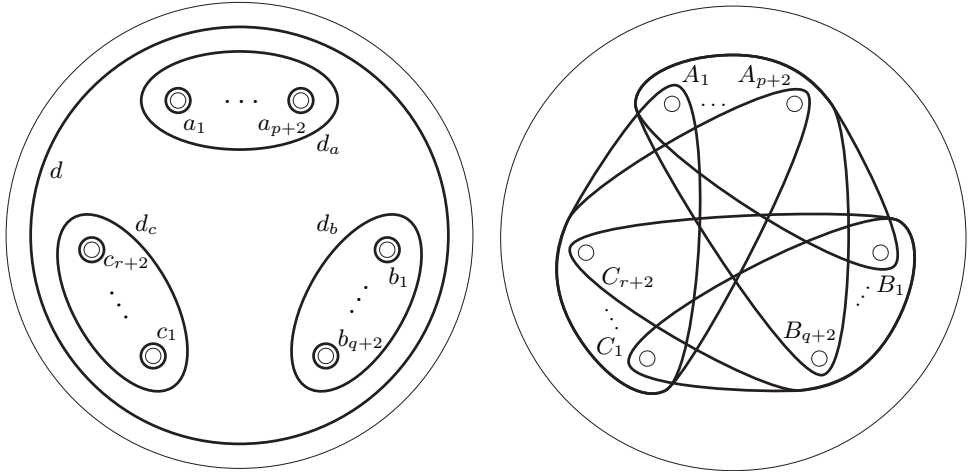
\begin{figure}[h]
\includegraphics[width=2.45in]{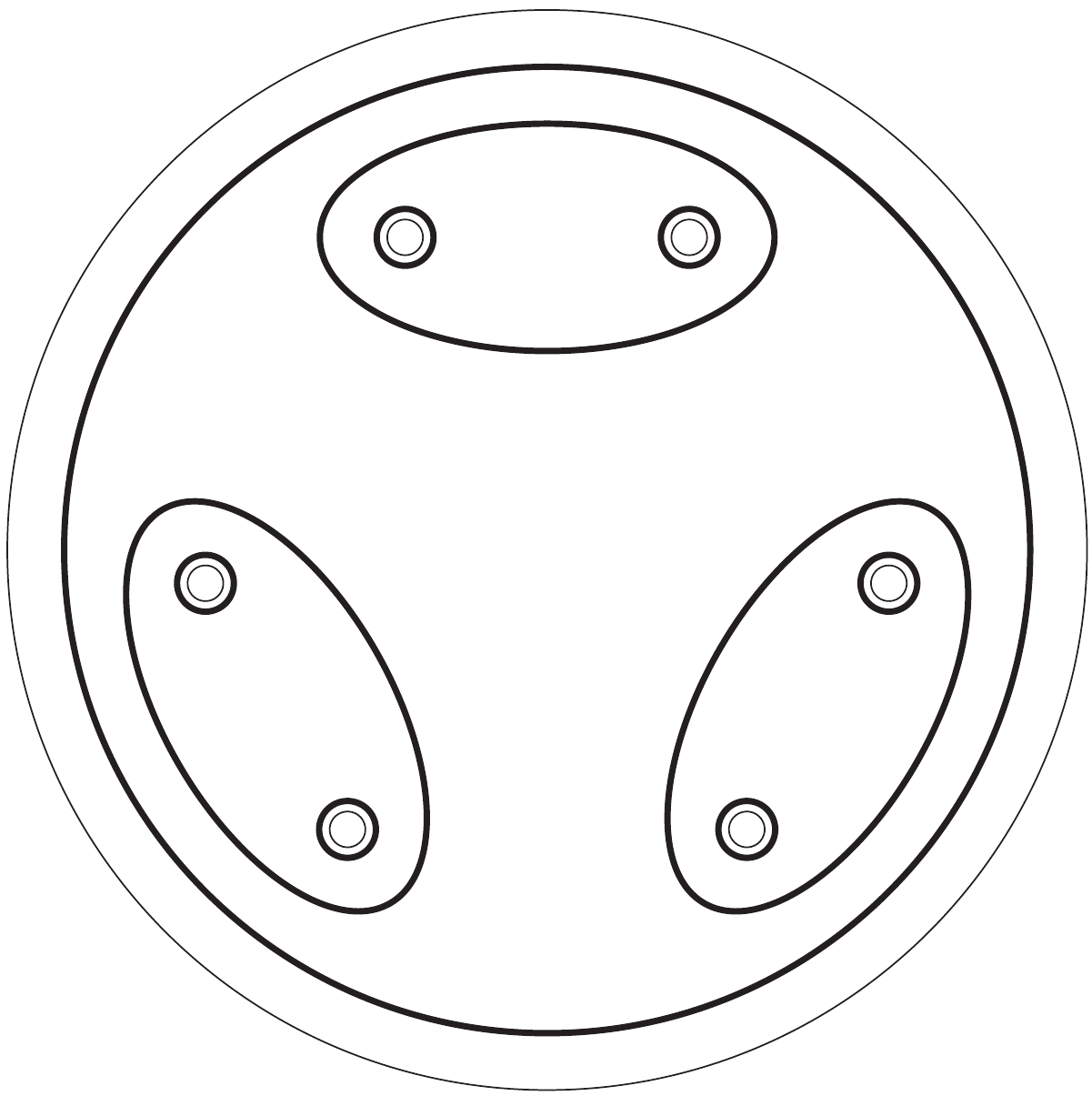}\quad
\includegraphics[width=2.45in]{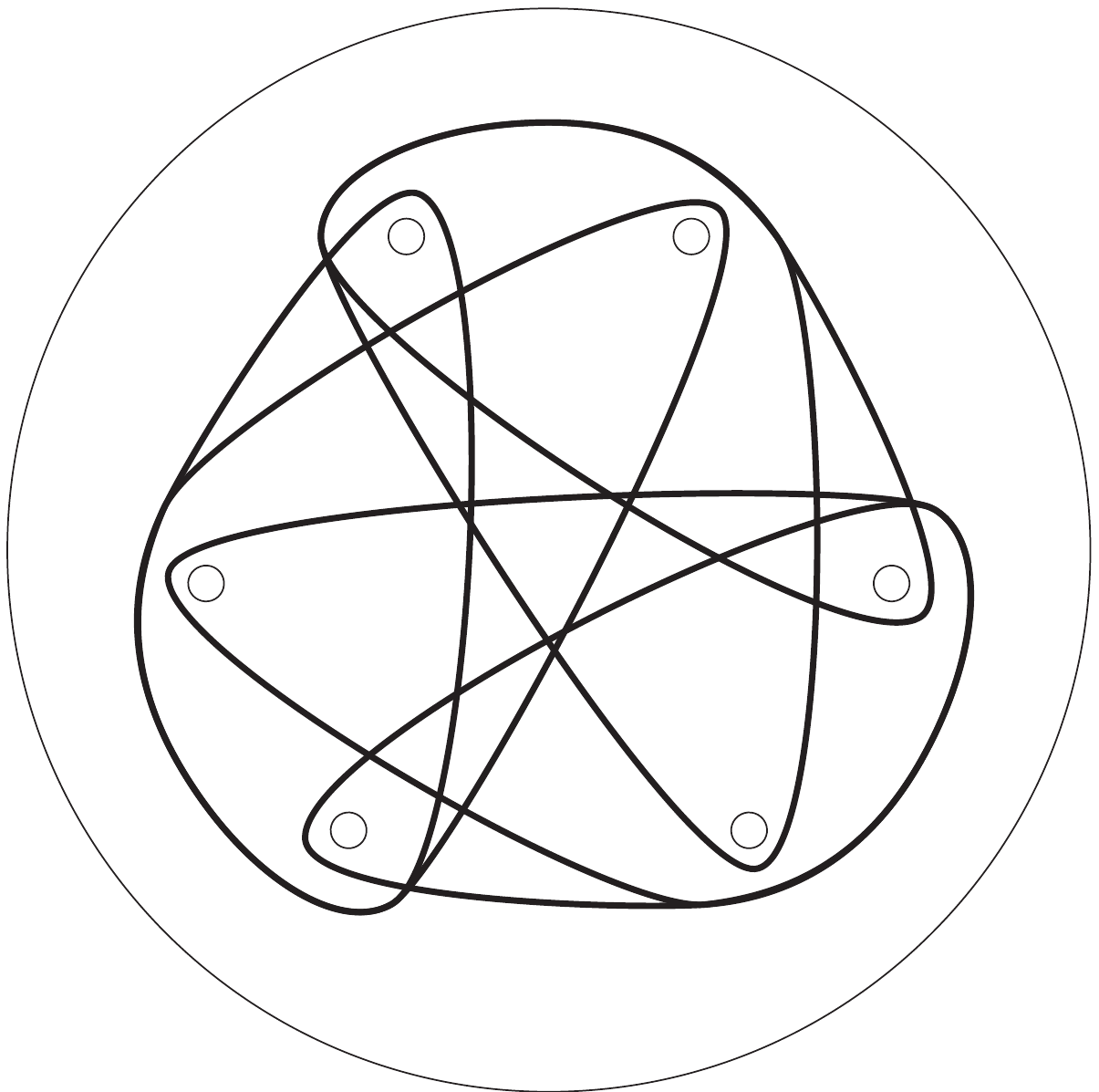}
\put(-295,128){\small $a_1$}
\put(-280,138){\large$\cdots$}
\put(-267,128){\small $a_{p+2}$}
\put(-245,120){\small $d_a$}
\put(-325,77){\small $c_{r+2}$}
\put(-305,50){\small $c_1$}
\put(-237,47){\small $b_{q+2}$}
\put(-218,70){\small $b_1$}
\put(-245,90){\small $d_b$}
\put(-313,90){\small $d_c$}
\put(-345,110){\small $d$}
\put(-232,58){\begin{rotate}{60}\large$\cdots$\end{rotate}}
\put(-323,68){\begin{rotate}{-60}\large$\cdots$\end{rotate}}
\put(-101,137){$\cdots$}
\put(-108,146){\small $A_1$}
\put(-87,146){\small $A_{p+2}$}
\put(-35,67){\small $B_1$}
\put(-70,50){\small $B_{q+2}$}
\put(-140,44){\small $C_1$}
\put(-138,70){\small $C_{r+2}$}
\put(-42,61){\begin{rotate}{60}$\cdots$\end{rotate}}
\put(-137,61){\begin{rotate}{-60}$\cdots$\end{rotate}}
\caption{\label{Wrelfig}Relation in the mapping class group of a planar surface corresponding to the $\mathcal W$ family of rational blowdowns.}
\end{figure}
\begin{proof} Begin with the lantern relation of Figure \ref{lanternfig}, and apply the Key Lemma to the pair of holes enclosed by $b$ and $c$ to split $c$. Since $b$ encloses just one hole of the surface, the commutation requirement of the Key Lemma is satisfied. We obtain a relation on a disk with four holes, of the form $ab^2c_1c_2d = xyC_2C_1$.
\begin{figure}[b]
\includegraphics[width=5in]{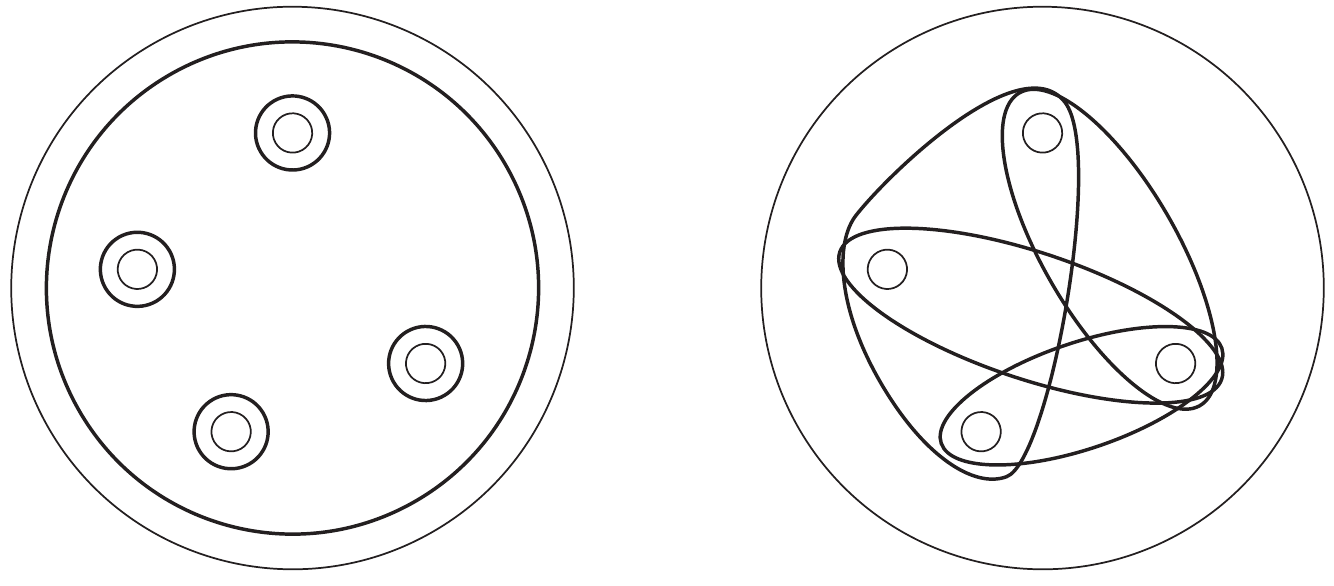}
\put(-318,69){\begin{rotate}{-60}\large$\cdots$\end{rotate}}
\put(-312,85){$c_{r+2}$}
\put(-295,55){$c_1$}
\put(-284,105){$a$}
\put(-264,58){$b$}
\put(-332,112){$d$}
\put(-112,64){\begin{rotate}{-60}\large$\cdots$\end{rotate}}
\put(-117,50){\small$C_1$}
\put(-112,70){\small$C_{r+2}$}
\put(-115,123){$x$}
\put(-55,123){$y$}
\put(-240,-5){$ab^{r+2}c_1\cdots c_{r+2}d = xyC_{r+2}\cdots C_1$}
\caption{\label{rel1}The lantern relation after $r+1$ applications of the Key Lemma.}
\end{figure}
 Applying the same procedure to the holes $b$ and $c_2$ gives $ab^2c_1c_2c_3d=xyC_3C_2C_1$ on a disk with five holes, and inductively we get the relation
\[
ab^{r+2}c_1\cdots c_{r+2}d = xyC_{r+2}\cdots C_1
\]
on a disk with $r+4$ holes as in Figure \ref{rel1} (this is, incidentally, the daisy relation of the introduction). 

Repeat this procedure using the holes $a$ and $b$ to split $b$ into $q+2$ holes and replace $y$ in the relation above by a product $B_{q+2}\cdots B_1$ around curves each enclosing $a$ and one of the holes $b_j$. We find 
\[
a^{q+2}b_1\cdots b_{q+2}c_1\cdots c_{r+2}d_b^{r+1}d = xB_{q+2}\cdots B_1C_{r+2}\cdots C_1,
\]
where $d_b$ is a circle enclosing all the $b_j$.

Finally we use the curve $x$ to split $a$ into $p+2$ holes. Observe that in this case one of the holes of the subsurface $S_3$ from the Key Lemma encloses more than one hole of our surface; however the commutation requirement is still trivial here since all the curves on the left hand side of the relation above are disjoint so that the corresponding twists pairwise commute, while the word $w_1'$ is trivial in this case. The desired relation follows.
\end{proof}

Clearly, the relation just obtained satisifies conditions (1), (2), and (3) mentioned previously.

\subsection{Relations for family $\mathcal N$}

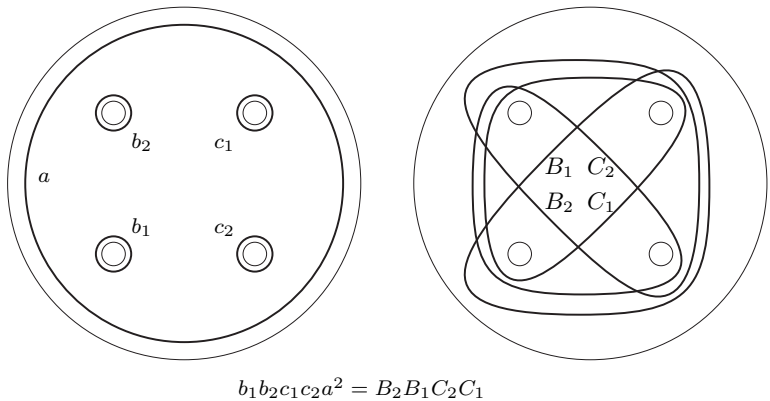
\begin{figure}[h]
\includegraphics[width=4in]{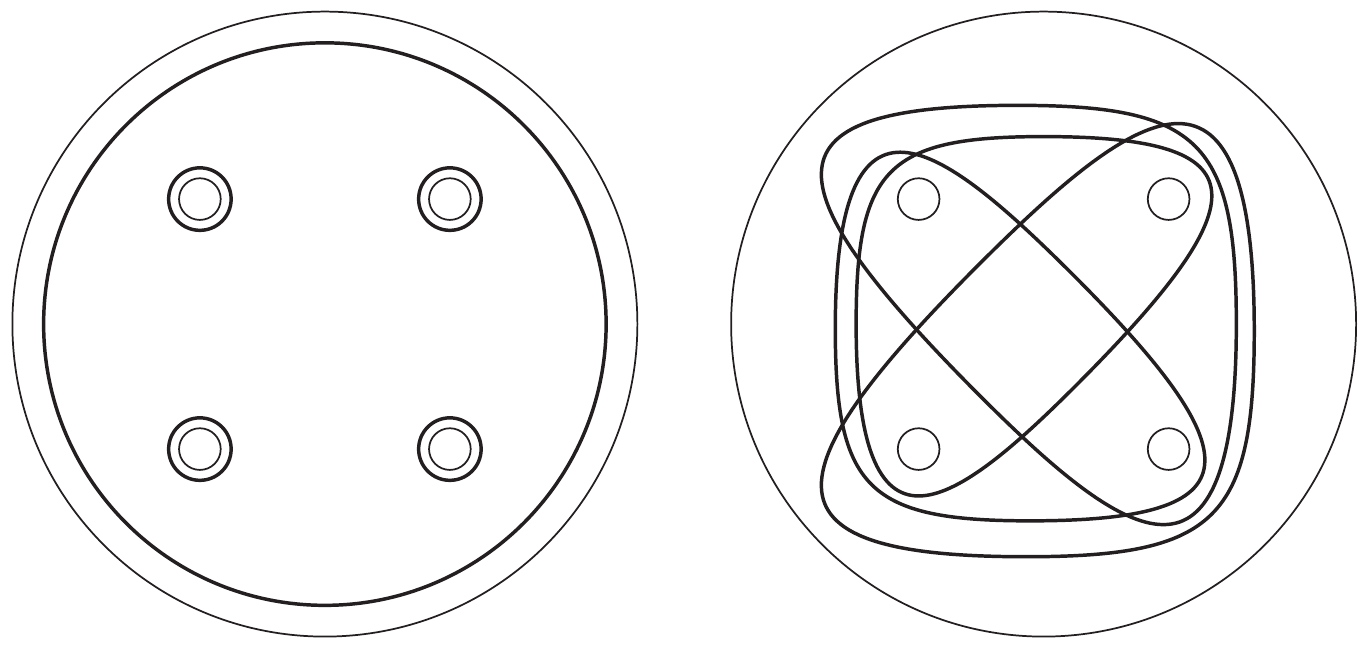}
\put(-275,70){$a$}
\put(-240,83){$b_2$}
\put(-240,51){$b_1$}
\put(-209,83){$c_1$}
\put(-209,51){$c_2$}
\put(-85,73){\small$B_1$}
\put(-69,73){\small$C_2$}
\put(-85,60){\small$B_2$}
\put(-69,60){\small$C_1$}
\put(-200,-10){$b_1b_2c_1c_2a^2=B_2B_1C_2C_1$}
\caption{\label{N1}A rearrangement and relabeling of the case $r = 1$ of Figure \ref{rel1}. The right side includes four curves, each enclosing three of the four holes in the disk.}
\end{figure}

 To derive this family of relations, take $r = 0$ in Figure \ref{rel1} and redraw it so that the hole enclosed by $b$ becomes the outer boundary: we have the relation $b_1b_2c_1c_2a^2=B_2B_1C_2C_1$ in the labeling conventions of Figure \ref{N1}. For $r\geq 0$, apply the Key Lemma consecutively $r+1$ times, using the highest-numbered curve $B_i$ to split the corresponding hole $b_i$ into two (the lemma applies since $w_1' = 1$ here). The result is the relation $b_1\cdots b_{r+3}c_1c_2d_c^{r+1}a^2 = B_{r+3}\cdots 
B_1 C_2C_1$, where $d_c$ is a curve enclosing the $c$-holes. 
 
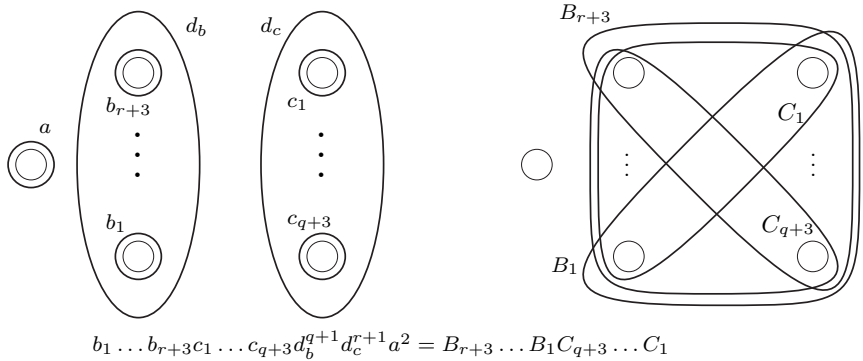
\begin{figure}[b]
\includegraphics[width=4.5in]{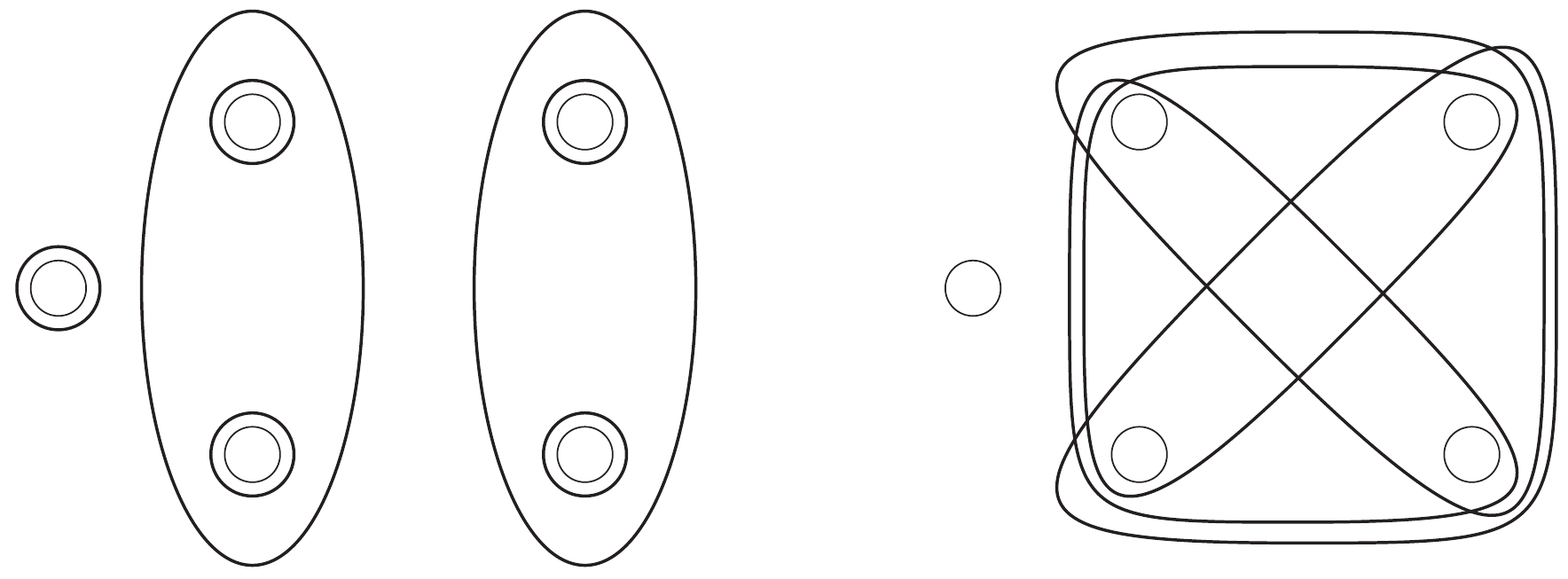}
\put(-310,72){$a$}
\put(-274,54){\begin{sideways}\LARGE$\cdots$\end{sideways}}
\put(-205,54){\begin{sideways}\LARGE$\cdots$\end{sideways}}
\put(-285,36){$b_1$}
\put(-285,80){$b_{r+3}$}
\put(-217,38){$c_{q+3}$}
\put(-217,82){$c_1$}
\put(-255,110){$d_b$}
\put(-227,110){$d_c$}
\put(-90,55){\begin{sideways}$\cdots$\end{sideways}}
\put(-20,55){\begin{sideways}$\cdots$\end{sideways}}
\put(-118,20){$B_1$}
\put(-115,115){$B_{r+3}$}
\put(-39,35){\small$C_{q+3}$}
\put(-33,77){\small$C_1$}
\put(-290,-10){$b_1\cdots b_{r+3}c_1\cdots c_{q+3}d_b^{q+1}d_c^{r+1}a^2 = B_{r+3}\cdots B_1C_{q+3}\cdots C_1$}
\caption{\label{N2}The case $p =0$ of a relation corresponding to the $\mathcal N$ family of blowdowns. (Drawn on two copies of a sphere with holes.) }
\end{figure}

Now for $q\geq 0$ repeat the process $q+1$ times using the highest-numbered $C_j$ to split the corresponding $c_1$ in two. Observe that the curve playing the role of $a$ in the Key Lemma is a circle enclosing all the holes $b_1,\ldots,b_{r+3}$, and the twist around this circle commutes with all the $C_j$. Hence the commutativity requirement of the Key Lemma holds. The result is a relation depicted in Figure \ref{N2}, which is drawn on a sphere with holes rather than a disk with holes. This will be the case $p=0$ of our family of relations corresponding to the family $\mathcal N$ of rational blowdowns.
 
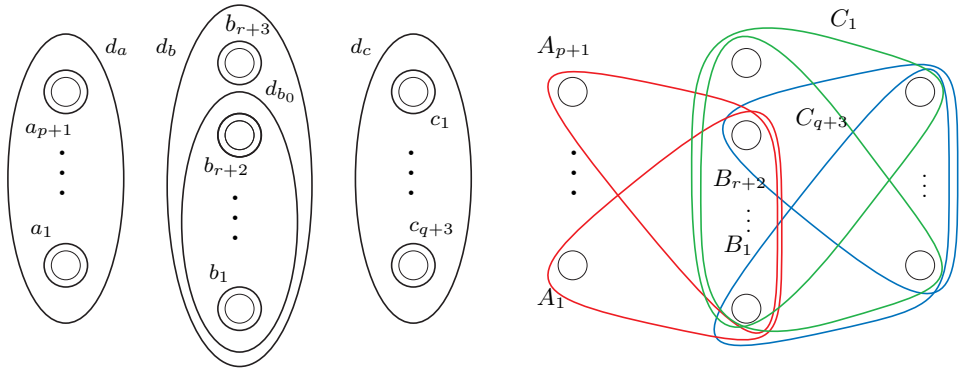
\begin{figure}[h]
\includegraphics[width=5in]{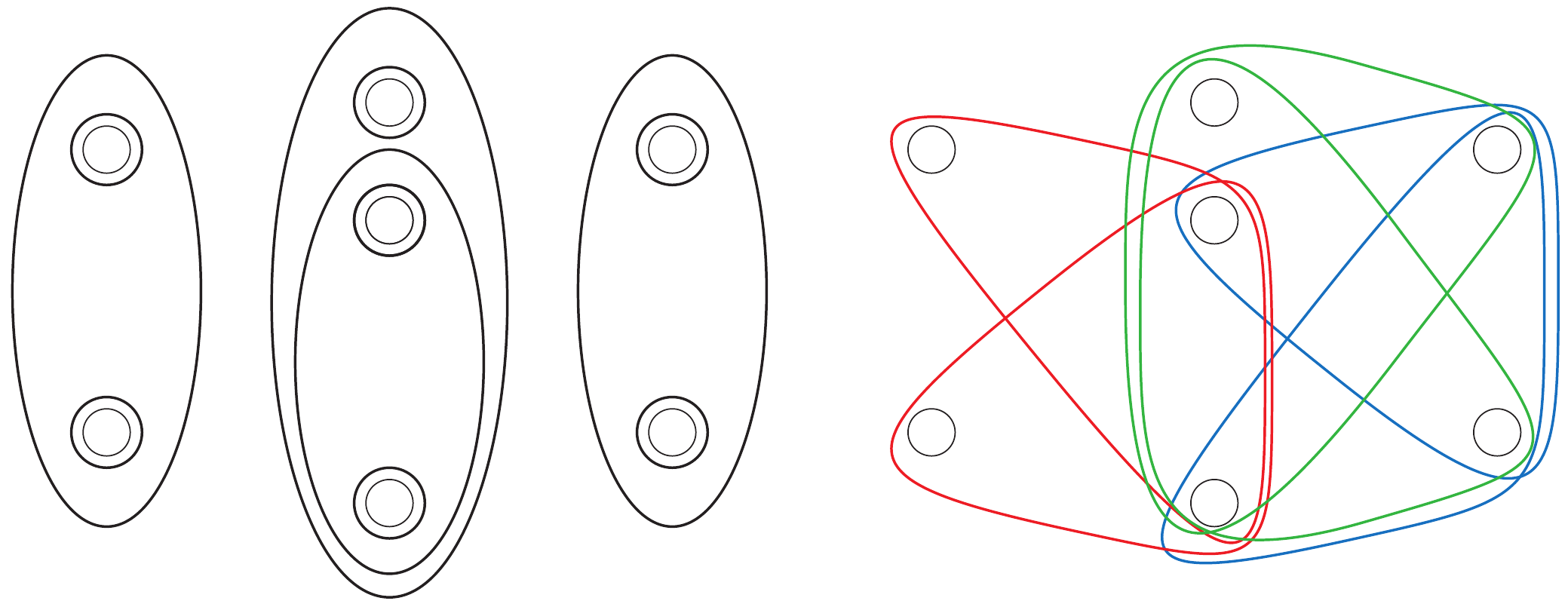}
\put(-322,120){$d_a$}
\put(-339,65){\begin{sideways}\LARGE$\cdots$\end{sideways}}
\put(-350,52){$a_1$}
\put(-352,90){$a_{p+1}$}
\put(-274,48){\begin{sideways}\LARGE$\cdots$\end{sideways}}
\put(-283,35){$b_1$}
\put(-285,75){$b_{r+2}$}
\put(-277,128){\small$b_{r+3}$}
\put(-262,104){$d_{b_0}$}
\put(-303,120){$d_b$}
\put(-208,65){\begin{sideways}\LARGE$\cdots$\end{sideways}}
\put(-208,53){$c_{q+3}$}
\put(-200,93){$c_1$}
\put(-230,120){$d_c$}
\put(-160,120){\small$A_{p+1}$}
\put(-147,65){\begin{sideways}\LARGE$\cdots$\end{sideways}}
\put(-160,25){\small$A_1$}
\put(-81,52){\begin{sideways}$\cdots$\end{sideways}}
\put(-90,45){\small$B_1$}
\put(-94,70){\small$B_{r+2}$}
\put(-63,93){\small$C_{q+3}$}
\put(-50,130){\small$C_1$}
\put(-15,65){\begin{sideways}$\cdots$\end{sideways}}
\caption{\label{N3}The relation for the general member of the $\mathcal N$ family. For clarity, curves $A_i$ are drawn in red, $B_j$ are blue, and $C_k$ are green.}
\end{figure}

Observe that the curve $B_{r+3}$ is isotopic (on the sphere with holes) to a curve $A$ enclosing the holes $b_1,\ldots, b_{r+2}$ and $a$. For given $p\geq 1$, apply the Key Lemma $p$ times to split hole $a$ into $a_1,\ldots, a_{p+1}$ (and $A$ into curves $A_1,\ldots A_{p+1}$), noting that here the element $w_1'$ from the Key Lemma is trivial: we obtain the relation of Figure \ref{N3}.
\begin{theorem}[$\mathcal N$-family of relations] Fix three integers $p,q,r\geq 0$, and let $F$ be a sphere with $p+q+r+7$ holes $a_1,\ldots, a_{p+1}$, $b_1,\ldots, b_{r+3}$, and $c_1,\ldots c_{q+3}$, and use the same letters to indicate boundary-parallel circles around these holes (and also Dehn twists around those curves). Let $d_b$ and $d_{b_0}$ be circles enclosing holes $b_1,\ldots, b_{r+3}$ and holes $b_1,\ldots,b_{r+2}$ respectively, and let $d_a$ and $d_c$ be circles enclosing all the holes labeled $a_i$ and all holes labeled $c_i$ respectively. Finally, let $A_i$, $B_j$, $C_k$ be circles such that:
\begin{enumerate}
\item For $i = 1,\ldots, p+1$, the curve $A_i$ encloses all holes $b_1,\ldots, b_{r+2}$ and hole $a_i$.
\item For $j=1,\ldots, r+3$, the curve $B_j$ encloses all holes $c_1,\ldots, c_{q+3}$ and hole $b_j$.
\item For $k = 1,\ldots, q+3$, the curve $C_k$ encloses holes $b_1,\ldots, b_{r+3}$ and hole $c_k$.
\end{enumerate}
See Figure \ref{N3}. Then if $p\geq 1$ we have the following relation in $\Mod(F,\partial F)$:
\[
a_1\cdots a_{p+1}b_1\cdots b_{r+3}c_1\cdots c_{q+3}d_ad_{b_0}^pd_b^{q+1}d_c^{r+1} = A_{p+1}\cdots A_1B_{r+2}\cdots B_1 C_{q+3}\cdots C_1.
\]
If $p=0$, then, writing $a$ for $a_1$,
\[
a^2b_1\cdots b_{r+3} c_1\cdots c_{q+3} d_b^{q+1}d_c^{r+1} = B_{r+3}\cdots B_1C_{q+3}\cdots C_1.
\]
\end{theorem}\hfill$\Box$

\subsection{Relations corresponding to linear plumbings}

Recall that given relatively prime integers $p > q>0$, there is a unique continued fraction expansion 
\[
-\frac{p^2}{pq-1} = [b_1,\ldots, b_k] = b_1 - \frac{1}{b_2 - \frac{1}{\cdots - \frac{1}{b_k}}}
\]
where each $b_i \leq -2$. Moreover, this expansion can be obtained from the base case $[-4]$ by repeated applications of the two operations 
\begin{enumerate}
\item[a)] $[b_1,\ldots,b_k] \mapsto [b_1-1,b_2,\ldots, b_k, -2]$
\item[b)] $[b_1,\ldots,b_k] \mapsto [-2, b_1,\ldots, b_{k-1}, b_k-1]$.
\end{enumerate} See, for example, \cite{SSW}, Proposition 4.1.
We indicate here how to obtain relations in the mapping class groups of planar surfaces by following a parallel procedure. The base case is the lantern relation itself. Supposing the construction of the continued fraction $[b_1,\ldots,b_k]$ to have begun with $n_1$ applications of the move (a) above, we apply the Key Lemma $n_1$ times to obtain the relation of Figure \ref{rel1}, with $r+2$ replaced by $n_1+1$. If the next stage calls for $n_2$ applications of (b), apply the Key Lemma $n_2$ times to the last-constructed curve in the previous stage, $C_{n_1+1}$, to split hole $b$ into $n_2+1$ holes.

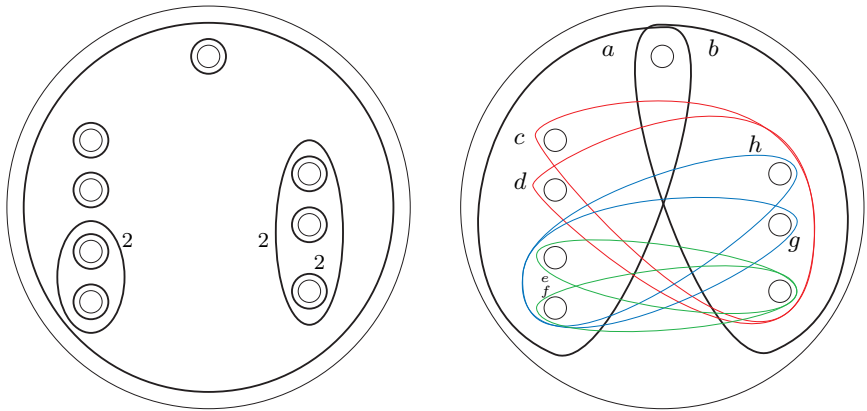
\begin{figure}[h]
\includegraphics[width=4.5in]{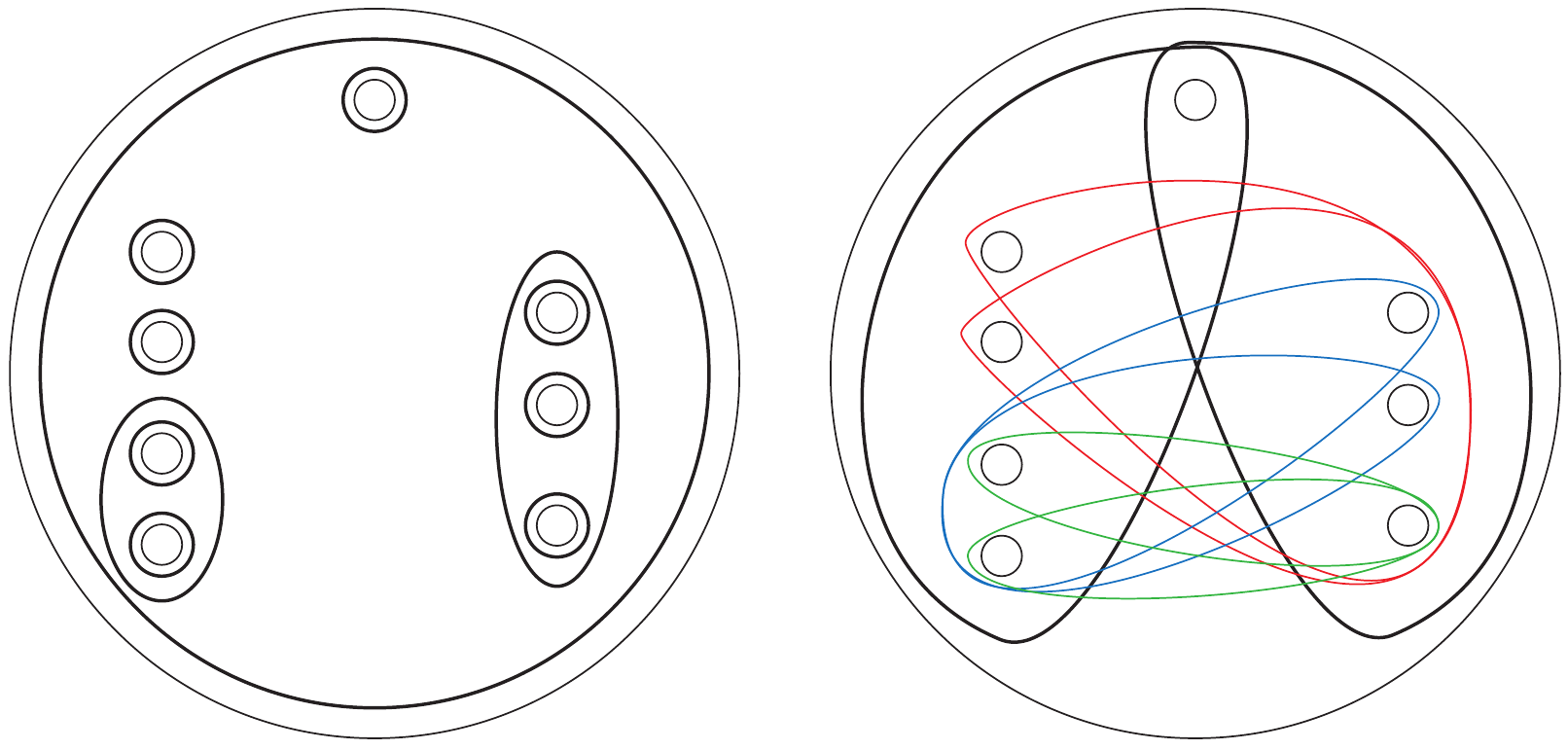}
\put(-280,63){$2$}
\put(-229,63){$2$}
\put(-208,55){$2$}
\put(-100,135){\small$a$}
\put(-60,135){\small$b$}
\put(-133,102){\small$c$}
\put(-133,85){\small$d$}
\put(-123,50){\tiny$e$}
\put(-123,45){\tiny$f$}
\put(-30,63){\small$g$}
\put(-45,99){\small$h$}
\caption{\label{examplefig}A relation corresponding to a rational blowdown along a linear plumbing diagram. The product of twists around the curves on the left (with the indicated curves taken more than once) is isotopic to the product of twists around the curves on the right, in aphabetical order. The colors indicate the stages in which the curves were constructed: the black are ``original'' to the lantern relation, then red, blue and green. Each application of the Key Lemma generates a new curve on the right diagram that we take to be the lower of the two possibilities.}
\end{figure}

And so on: the next stage of construction of the continued fraction calls for, say, $n_3$ applications of (a). The parallel is to apply the Key Lemma $n_3$ times using the last-constructed curve to split the hole that was not split during the previous stage. The process continues in this manner to obtain a relation in a planar mapping class group corresponding to each fraction $p^2/(pq-1)$. As an example, the relation corresponding to $p=17$ and $q= 7$, with continued fraction $[-3,-2,-6,-2,-4,-2]$, is shown in Figure \ref{examplefig}. Note that here $n_1 = n_2 = 2$, $n_3 = 1$.


\section{Kirby diagrams for rational balls} 
\label{diagramsec}

The proof that the monodromy substitutions outlined previously result in rational blowdowns amounts to analyzing Kirby pictures of Lefschetz fibrations with bordered fibers over $D^2$: specifically, the Lefschetz fibrations with planar pages arising in our monodromy relations. In the following, we consider Kirby pictures for the 4-manifolds described by each side of the various relations and show that one side describes a plumbing of spheres while the other side is a rational ball. (The second of these points is essentially obvious, but we give explicit Kirby pictures in any case to relate these rational blowdowns to those in \cite{FS} and \cite{SSW}.)

To begin with, we observe that if $F$ is a planar surface, diffeomorphic to a disk with $n$ smaller disks removed, and if $W$ is a Lefschetz fibration over $D^2$ having standard fiber $F$ and described by the global monodromy word $w$, then there is a standard Kirby diagram for $W$ obtained as follows. First note that a neighborhood $D^2\times F$ of a regular fiber is diffeomorphic to the complement of $n$ standard disks in $D^2 \times D^2 = B^4$, and therefore has a Kirby diagram consisting of an $n$-component unlink decorated by dots. In fact, drawing $F$ as a round disk in the $xy$ plane with $n$ disks deleted, we can consider the components of this unlink as vertical line segments passing through the centers of the $n$ deleted disks (and closed in a standard way). Indeed, the decomposition of $B^4$ into $D^2\times D^2$ corresponds to the usual decomposition of $S^3$ into the union of solid tori, a meridian disk of one of which we identify with the disk containing $F$. Vertically-displaced parallel copies of $F$ apparent in the picture are then nearby fibers. The standard disks we must remove from $D^2\times D^2$ then have boundary equal to the circles swept out by the centers of the disks to delete from $F$, as we move around one of the solid tori. 

To build $W$ from this point, we add $-1$ framed 2-handles to vanishing circles on parallel copies of $F$, in the order (from bottom to top) in which they appear in the word $w$. In particular, if the collection of vanishing cycles has the property that their rational homology classes form a basis for $H_1(F;\cue)$, then it is easy to see that $W$ has trivial rational homology. Furthermore, the effect of a monodromy substitution is local: it amounts to replacing a piece of a global Lefschetz fibration $X(w_1\cdot w')$ diffeomorphic to $W(w_1)$ by another piece $W(w_2)$ (here we think of $W(w_i)$ as Lefschetz fibrations whose fiber is a subsurface of that of $X$---a planar subsurface in all the cases at hand). Thus we need only understand the topology of $W(w_i)$: the properties (1), (2), and (3) of the previous section ensure that $W(w_1)$ is diffeomorphic to a plumbing (but see below for a direct verification of this), and that $W(w_2)$ is a rational ball.

 \subsection{Graphs in family $\W$} 
\begin{figure}[h]
\includegraphics[width=5in]{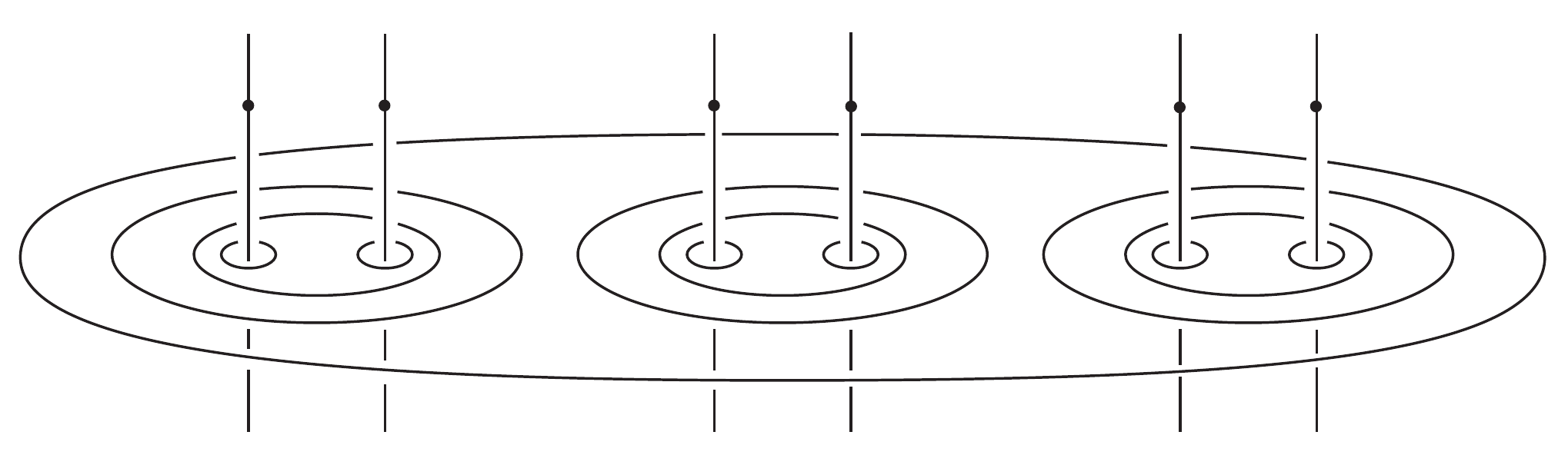}
\put(-306,105){$\overbrace{\hspace{.5in}}^{p+2}$}
\put(-294,48){$\cdots$}
\put(-257,48){$\cdots$}
\put(-258.5,53.5){\tiny $q\hspace{-.8ex}+\hspace{-.8ex}1$}
\put(-199,105){$\overbrace{\hspace{.5in}}^{q+2}$}
\put(-187,48){$\cdots$}
\put(-150,48){$\cdots$}
\put(-151.5,53.5){\tiny $r\hspace{-.8ex}+\hspace{-.8ex}1$}
\put(-92,105){$\overbrace{\hspace{.5in}}^{r+2}$}
\put(-80,48){$\cdots$}
\put(-43,48){$\cdots$}
\put(-44.5,53.5){\tiny $p\hspace{-.8ex}+\hspace{-.8ex}1$}
\caption{\label{figAA}Handle diagram for the Lefschetz fibration described by one side of the $\mathcal W$-relation. Unlabeled circles have framing $-1$.}
\end{figure}
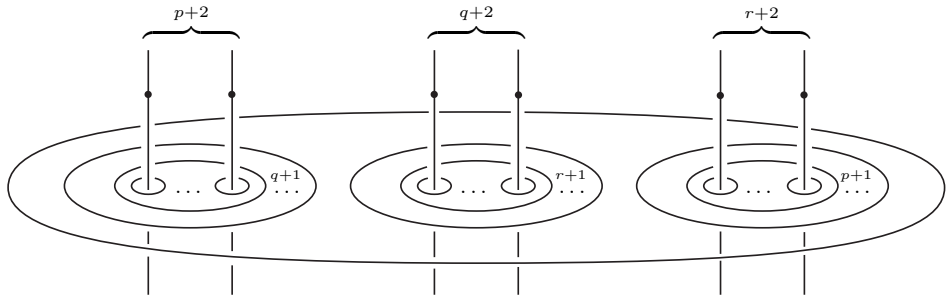

\begin{figure}[h]
\includegraphics[width=5in]{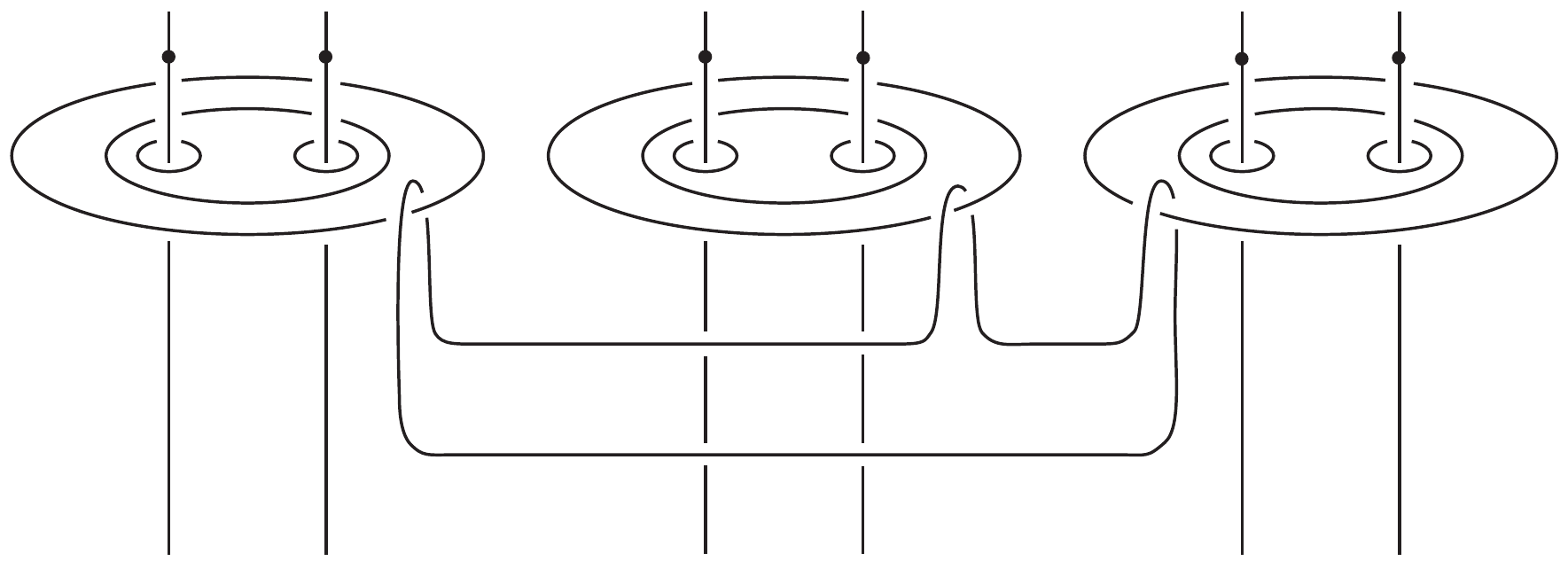}
\put(-314,123){\huge$\cdots$}
\put(-311,136){${}_{p+2}$}
\put(-190,123){\huge$\cdots$}
\put(-187,136){${}_{q+2}$}
\put(-66,123){\huge$\cdots$}
\put(-63,136){${}_{r+2}$}
\put(-353,93){$\cdots$}
\put(-353,86){${}^{q+1}$}
\put(-310,93){$\cdots$}
\put(-267,93){$\cdots$}
\put(-229,93){$\cdots$}
\put(-229,86){${}^{r+1}$}
\put(-186,93){$\cdots$}
\put(-143,93){$\cdots$}
\put(-106,93){$\cdots$}
\put(-63,93){$\cdots$}
\put(-20,93){$\cdots$}
\put(-21,86){${}^{p+1}$}
\put(-240,15){$-4$}
\caption{\label{AB}Unlabeled 2-handles carry framing $-1$.}
\end{figure}
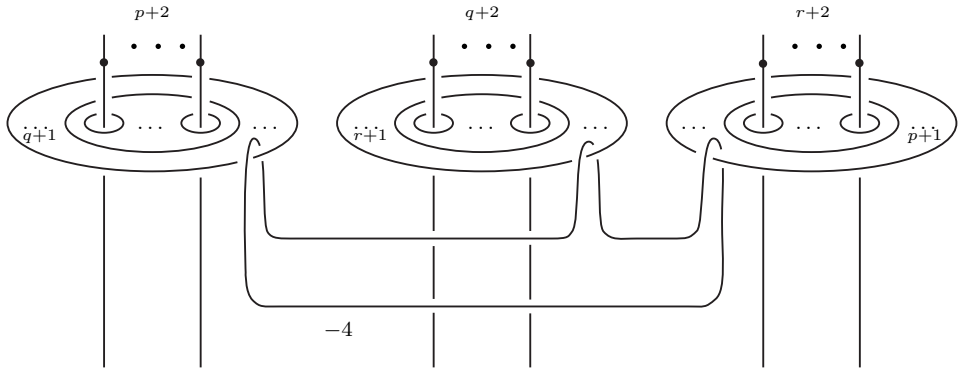

\begin{figure}[h]
\includegraphics[width=4.4in]{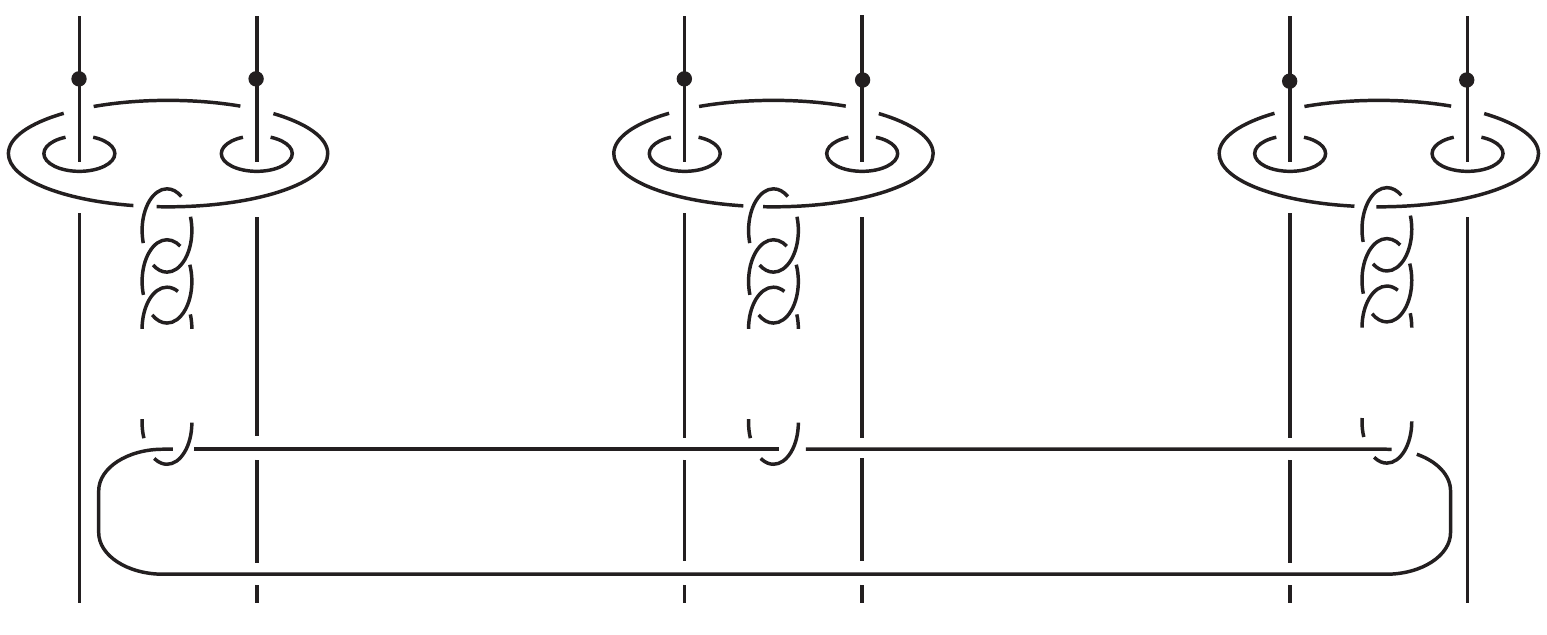}
\put(-294,121){\huge$\cdots$}
\put(-291,134){${}_{p+2}$}
\put(-169,121){\huge$\cdots$}
\put(-166,134){${}_{q+2}$}
\put(-44,121){\huge$\cdots$}
\put(-41,134){${}_{r+2}$}
\put(-290,97){$\cdots$}
\put(-165,97){$\cdots$}
\put(-40,97){$\cdots$}
\put(-285,49){$\vdots$}
\put(-279,51){$q$}
\put(-160,49){$\vdots$}
\put(-154,51){$r$}
\put(-34,49){$\vdots$}
\put(-28,51){$p$}
\put(-230,43){$-4$}
\caption{\label{AC}The small circles in the vertical chains carry framing $-2$, while the remaining unlabeled 2-handles carry framing $-1$.}
\end{figure}
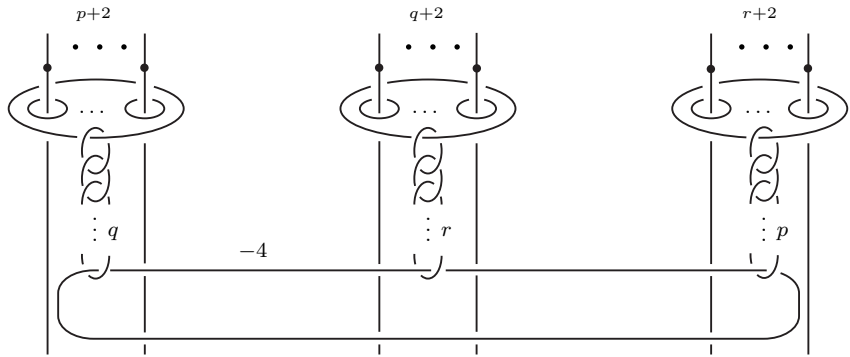

We first verify that the Lefschetz fibration described by the longer of our two words in the relation for the $\W$-family of substitutions is diffeomorphic to a plumbing of disk bundles over spheres according to the graph $\Gamma_{p,q,r}$ of \cite{SSW}. The procedure outlined above gives rise to the Kirby picture in Figure \ref{figAA} (note that all the twists in this word commute, so the order of the handle additions is immaterial).
Slide the outermost $-1$ circle over the three 2-handles that are outermost in their respective families to obtain Figure \ref{AB}.
 Then successively slide each of the outermost 2-handles in the three families over the next one in to get a chain of $-2$-framed circles as shown in Figure \ref{AC}.
 Now slide the last remaining ``large'' 2-handles over each remaining $-1$ circle in its family, cancelling the latter with its corresponding 1-handle at each stage. This has the effect of subtracting the number of 1-handles in that cluster from the framing on the last 2-handle; the end result is the desired plumbing.
 
 Turning to the word on the right hand side of the relation of Theorem \ref{Wrelthm}, note first that since the length of this word is equal to $b_1(F)$ and since the vanishing cycles are easily seen to span $H_1(F;\cue)$, the fact that the corresponding Lefschetz fibration over $D^2$ is a rational ball is clear. For the sake of completeness, however, we derive an explicit picture of this rational ball, and show it is diffeomorphic to one constructed in \cite{SSW}. In fact, we begin with the latter construction. 
 
\begin{figure}[b]
\includegraphics[width=3in]{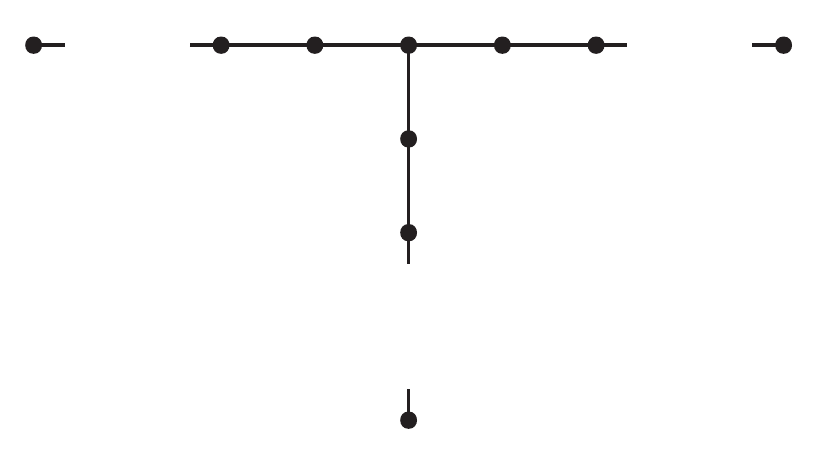}
\put(-212,110){$\underbrace{\hspace{.8in}}_{p+1}$}
\put(-189,112.5){$\cdots$}
\put(-150,119){${}^{-(q+2)}$}
\put(-110,122){\scriptsize $1$}
\put(-98,119){${}^{-(r+2)}$}
\put(-62,110){$\underbrace{\hspace{.8in}}_{q+1}$}
\put(-39,112.5){$\cdots$}
\put(-141,85){${}^{-(p+2)}$}
\put(-109.5,35){$\vdots$}
\put(-115,37){$\left.\begin{array}{c}\vspace*{.68in} \\ \end{array}\right\}$}
\put(-95,34){${}^{r+1}$}
\caption{\label{AJ}}
\end{figure}

 Recall from \cite{SSW} that to a negative-definite star-shaped plumbing graph $\Gamma$ one can associate a ``dual graph'' $\Gamma'$ such that the corresponding plumbed 4-manifolds $W(\Gamma)$ and $W(\Gamma')$  have diffeomorphic boundaries $Y(\Gamma)$ with opposite orientation, and $W(\Gamma)\cup_{Y(\Gamma)} W(\Gamma')$ is diffeomorphic to a blowup of $\cee P^2$. In certain cases one can also find the dual plumbing embedded in a (different) blowup of $\cee P^2$ in such a way that $\Gamma'$ spans the rational homology; its complement in this embedding is then a rational ball bounded by $Y(\Gamma)$ with the correct orientation. This plan is carried out in \cite{SSW}; in that case the dual graph to $\Gamma_{p,q,r}$ is a star-shaped graph with three legs shown in Figure \ref{AJ}, where the undecorated vertices have weight $-2$.
 Following \cite{SSW} we find such a plumbing in a rational surface by blowing up a configuration of four general lines in $\cee P^2$.

\begin{figure}[h]
\includegraphics[width=3in]{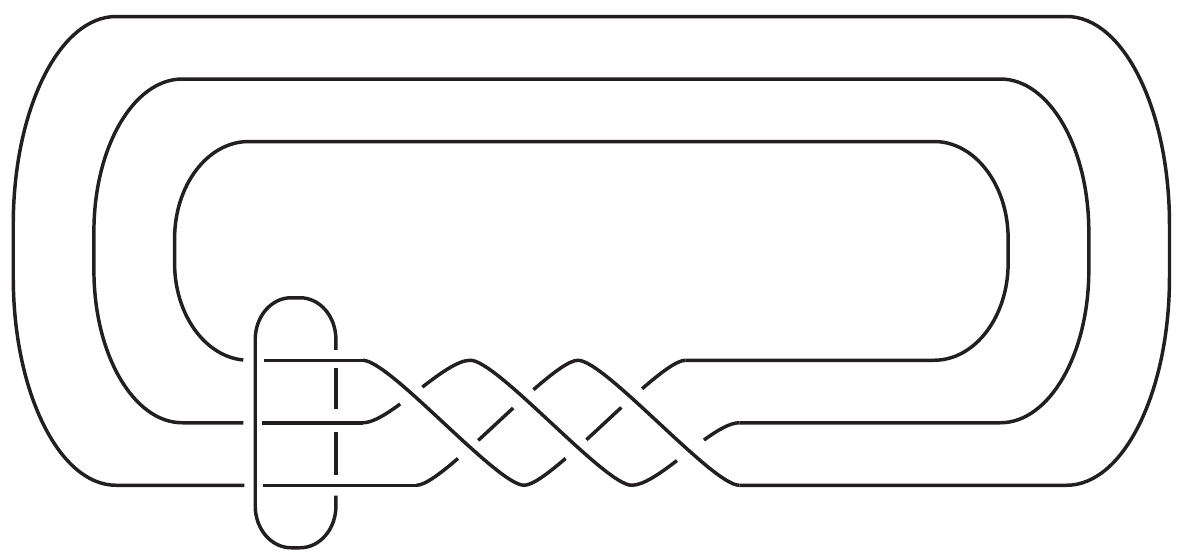}
\put(-228,53){$L_4$}
\put(-212,53){$L_3$}
\put(-197,53){$L_2$}
\put(-170,53){$L_1$}
\put(-5,5){\begin{tabular}{r} $\cup$\hspace{.5ex}  3 $3$-handles\\  4-handle\end{tabular}}
\caption{\label{AE}A Kirby picture for $\cee P^2$. All framings are $+1$.}
\end{figure}
\begin{figure}[t]
\includegraphics[width=4in]{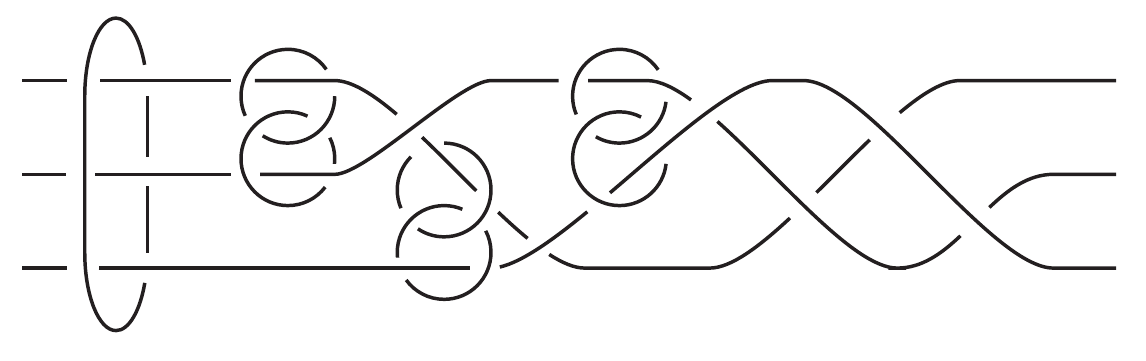}
\put(-30,60){$-2$}
\put(-30,36){$-2$}
\put(-30,12){$-2$}
\put(-260,95){$1$}
\put(-225,83){$-2$}
\put(-225,30){$-1$}
\put(-205,38){$*$}
\put(-170,57){$-1$}
\put(-163,42){$*$}
\put(-183,5){$-2$}
\put(-142,83){$-2$}
\put(-140,30){$-1$}
\put(-118,47){$*$}
\put(-35,-15){\begin{tabular}{r} $\cup$\hspace{.5ex}  3 $3$-handles\\  4-handle\end{tabular}}
\caption{\label{AF}The braided region of Figure \ref{AE} after six blowups. Additional blowups may be performed at the indicated crossings.}
\end{figure}

 Figure \ref{AE} is a Kirby picture for $\cee P^2$ in which four general lines are visible. We blow up twice each at the intersections between $L_2$, $L_3$, and $L_4$ to give Figure \ref{AF}. 
 Here the basic dual plumbing is visible (for $p=q=r=0$); $p$, $q$, and $r$ additional blowups at the indicated intersections allow us to construct the dual of any desired $\Gamma_{p,q,r}$. The complement of $\Gamma'$ is given by the three $-1$-framed 2-handles together with the 3- and 4-handles in Figure \ref{AF}. To obtain a standard Kirby picture for it, we must turn these handles upside down: this is complicated by the presence of 3-handles, which become 1-handles upon inverting the handlebody and should therefore be attached before the (inverted) 2-handles. 
 
 To obtain the picture, recall that the 2-handle dual to one of the $-1$ framed circles is attached along a 0-framed meridian $m$ to that circle. Now, the 3-manifold described by all the 2-handles in Figure \ref{AF} is necessarily diffeomorphic to $\# 3 S^1\times S^2$, and can therefore be simplified to a diagram of a 0-framed 3-component unlink. Replacing 0's by dots then describes the duals to the 3-handles. Our task is to follow the three meridians $m$ through this simplification; the resulting framed link together with the dotted unknots gives the desired diagram for our rational ball, after reversing its orientation.
 
\begin{figure}
\includegraphics[width=4.25in]{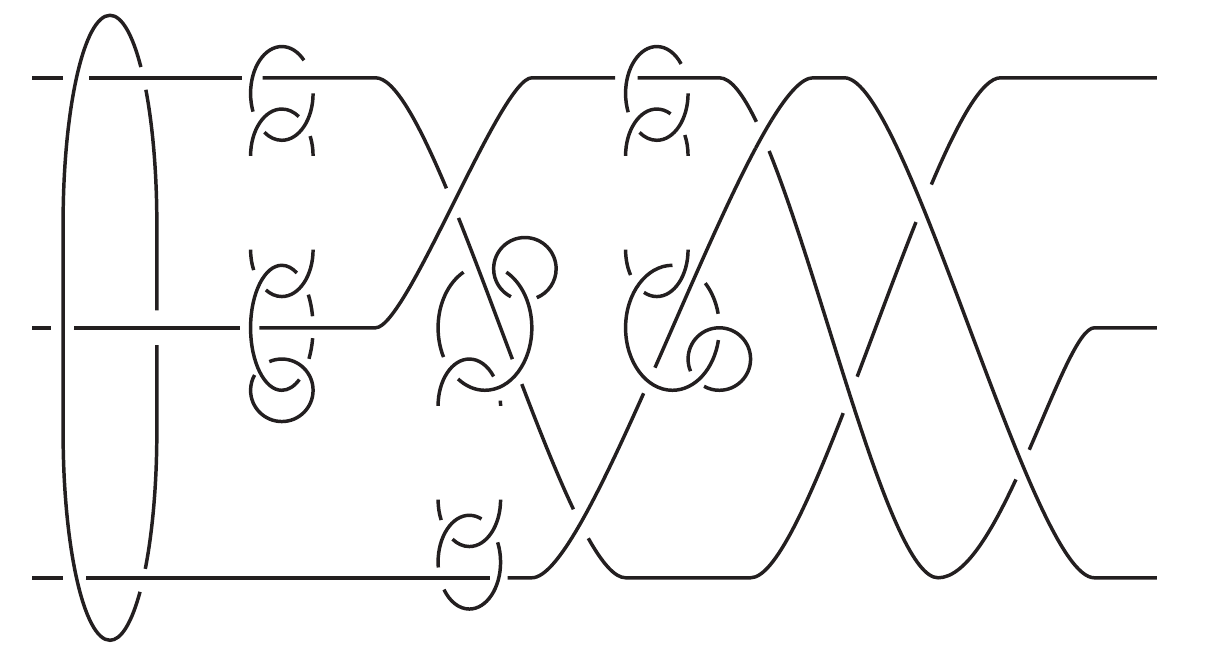}
\put(-50,153){$-(q+2)$}
\put(-50,93){$-(p+2)$}
\put(-50,10){$-(r+2)$}
\put(-280,-5){$1$}
\put(-237,112){$\vdots$}
\put(-260,90){$-1$}
\put(-245,50){$0$}
\put(-178,112){$0$}
\put(-170,85){$-1$}
\put(-190,48){$\vdots$}
\put(-143,112){$\vdots$}
\put(-163,70){$-1$}
\put(-120,61){$0$}
\caption{\label{AGa} Finding the plumbing complement. Unlabeled circles carry framing $-2$; there are $p+1$ such circles in the left column, $q+1$ in the center and $r+1$ on the right.}
\end{figure}

\begin{figure}
\includegraphics[width=4.25in]{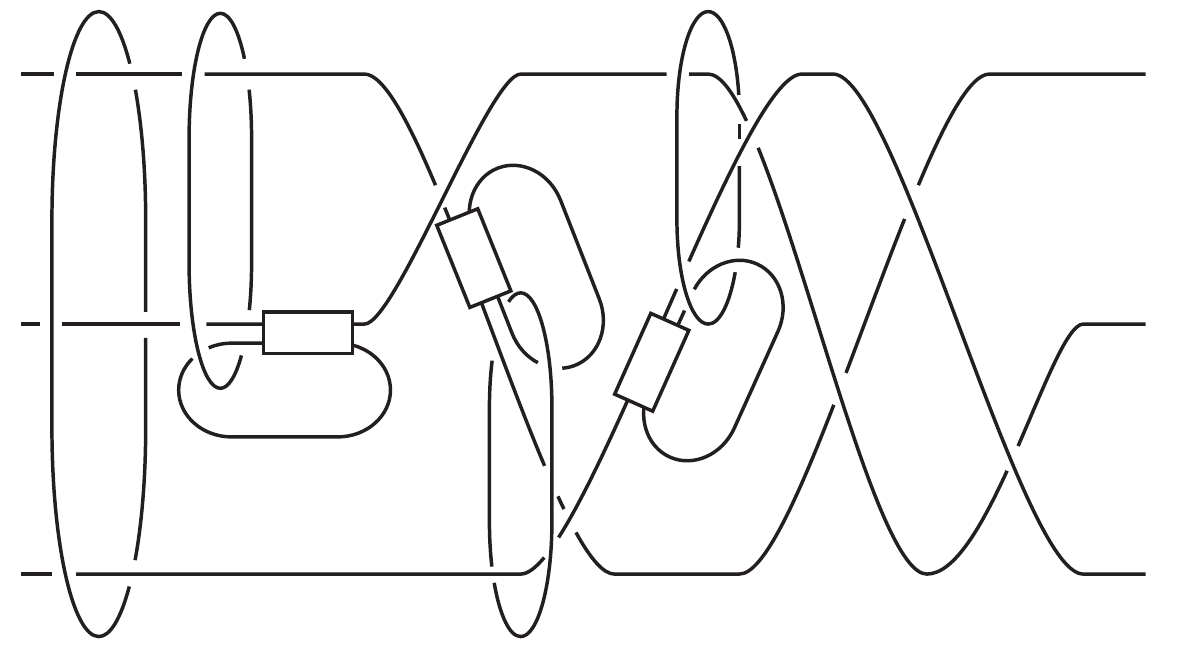}
\put(-236,80){${}^{p+1}$}
\put(-193,109){\begin{rotate}{-70}${}^{q+1}$\end{rotate}}
\put(-138.5,69){\begin{rotate}{66}${}^{r+1}$\end{rotate}}
\put(-30,155){$-1$}
\put(-30,93){$-1$}
\put(-30,10){$-1$}
\put(-280,-5){$1$}
\put(-260,172){$-1$}
\put(-245,50){$p+1$}
\put(-170,133){$q+1$}
\put(-185,-5){$-1$}
\put(-135,172){$-1$}
\put(-140,45){$r+1$}
\caption{\label{AGb}After blowing down most of the circles in the chains in Figure \ref{AGa}.}
\end{figure}

\begin{figure}
\includegraphics[width=4.25in]{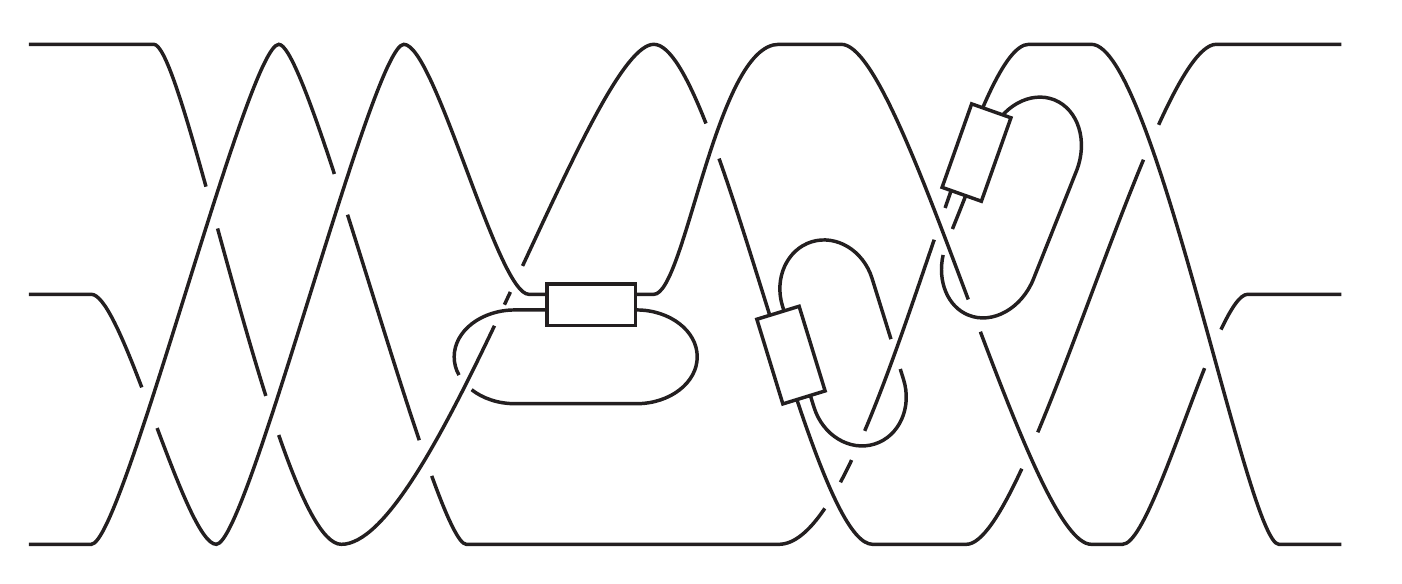}
\put(-185,56){${}^{p+2}$}
\put(-141,56){\begin{rotate}{-74}${}^{q+2}$\end{rotate}}
\put(-90.5,86){\begin{rotate}{70}${}^{r+2}$\end{rotate}}
\put(-190,30){$p+2$}
\put(-135,80){$q+2$}
\put(-83,106){${}^{r+2}$}
\put(-20,123){$0$}
\put(-20,70){$0$}
\put(-20,15){$0$}
\caption{\label{AGc}Blowing down remaining $\pm 1$ circles in Figure \ref{AGb}.}
\end{figure}

\begin{figure}
\includegraphics[width=3in]{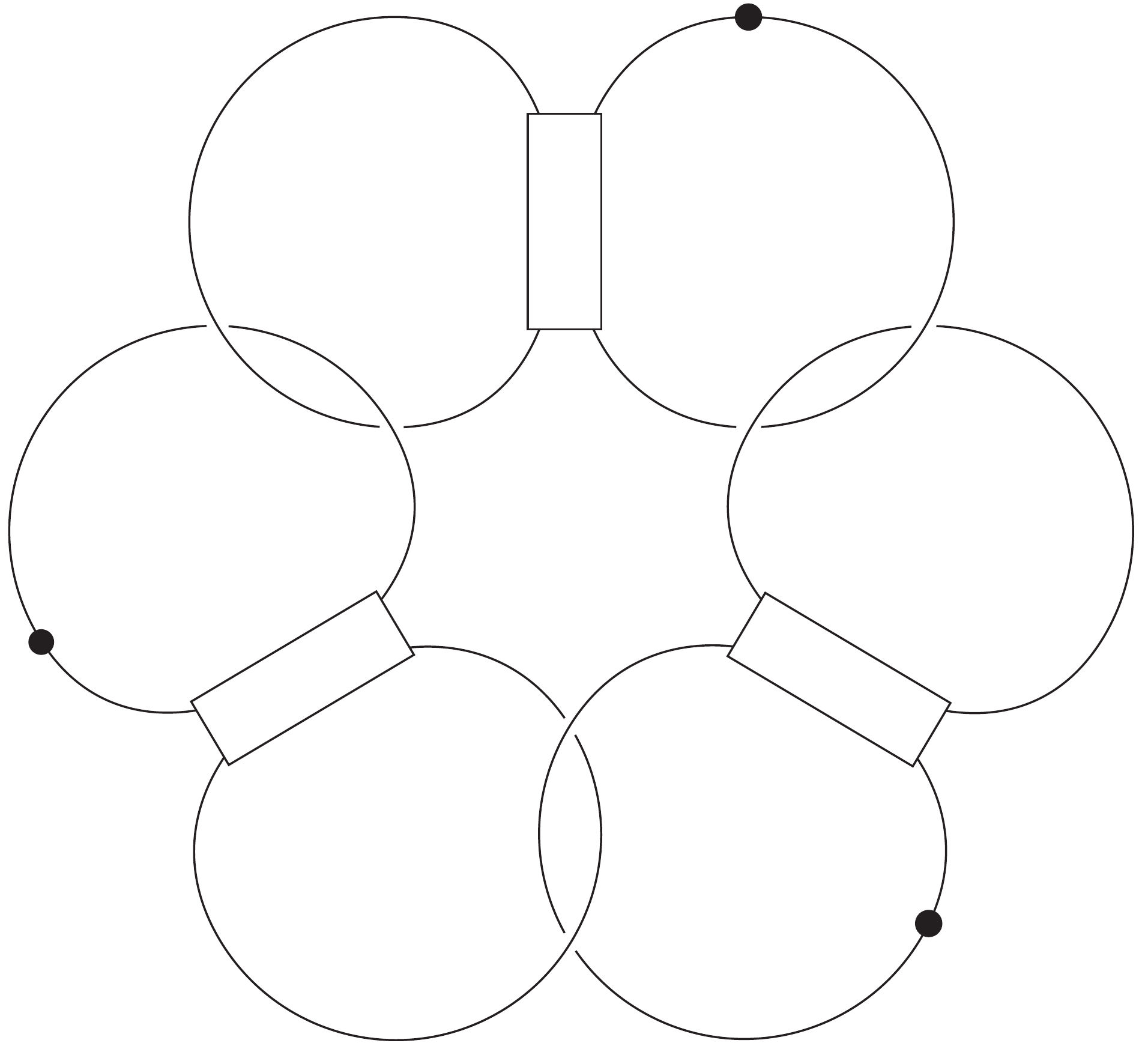}
\put(-114,145){\begin{sideways}${}^{-(p+2)}$\end{sideways}}
\put(-170,60){\begin{rotate}{31}${}^{-(q+2)}$\end{rotate}}
\put(-73,75){\begin{rotate}{-31}${}^{-(r+2)}$\end{rotate}}
\put(-190,200){$-(p+2)$}
\put(-190,-5){$-(q+2)$}
\put(0,105){$-(r+2)$}
\caption{\label{AGd}A rational ball bounded by the plumbed 3-manifold described by the graph $\Gamma_{p,q,r}$. This is obtained from Figure \ref{AGc} by closing the ends of the $0$-circles and exchanging $0$ framings for dots, followed by isotopy and reversing orientation.}
\end{figure}

 The sequence of Kirby moves is shown in Figures \ref{AGa}--\ref{AGd}, where we have passed to the case of general $p,q,r$, and performed the necessary orientation reversal only at the last step, Figure \ref{AGd}.

\begin{figure}[t]
\includegraphics[width=5in]{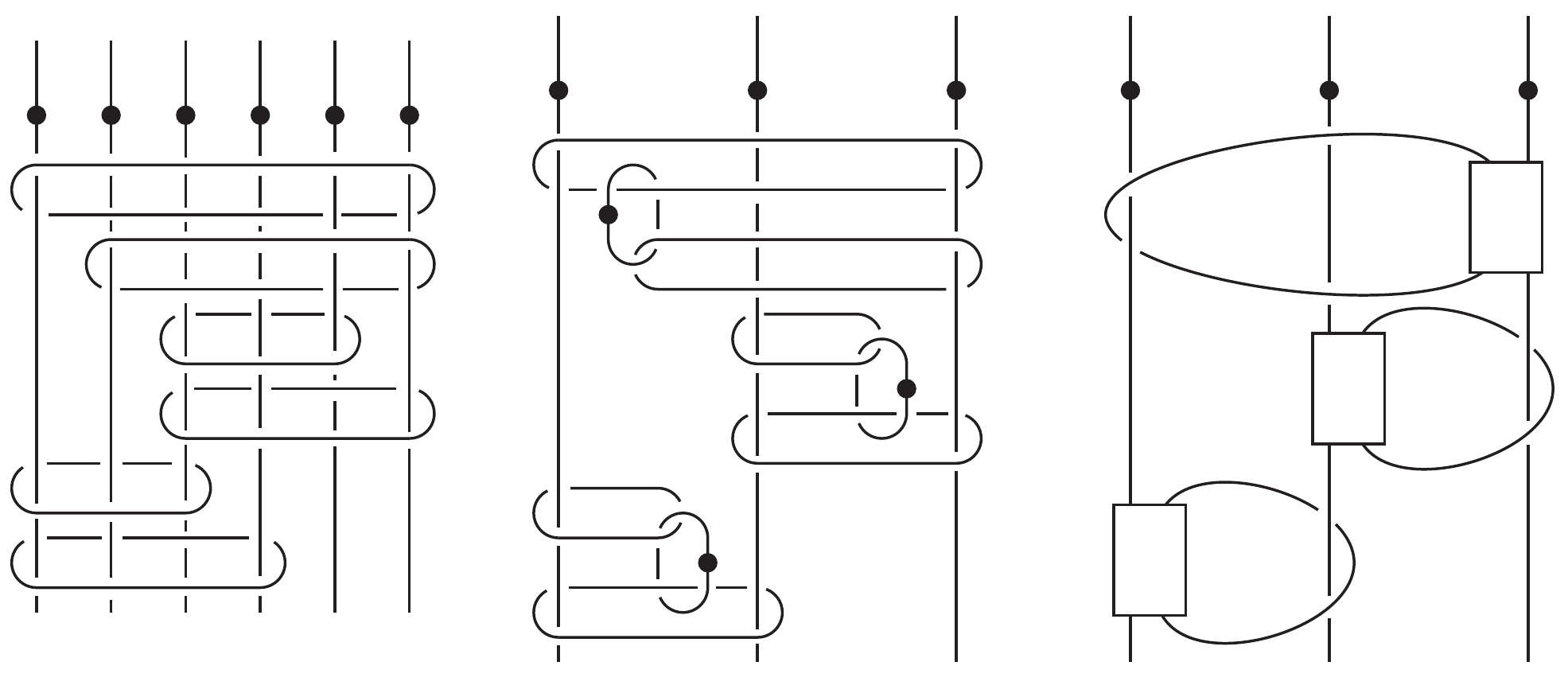}
\put(-316,-10){(a)}
\put(-195,-10){(b)}
\put(-60,-10){(c)}
\put(-352,138){${\longleftarrow}$}
\put(-318,138){$\longrightarrow$}
\put(-284,138){$\longrightarrow$}
\put(-200,108){$\uparrow$}
\put(-175,68){$\downarrow$}
\put(-220,27.5){$\downarrow$}
\put(-95,83){${}^{-(p+2)}$}
\put(-95,50){${}^{-(q+2)}$}
\put(-45,40){${}^{-(r+2)}$}
\put(-100,18){\begin{sideways} \small${}^{-(q+2)}$\end{sideways}}
\put(-54,57.5){\begin{sideways} \small${}^{-(r+2)}$\end{sideways}}
\put(-17,97){\begin{sideways} \small${}^{-(p+2)}$\end{sideways}}
\caption{\label{AI}Simplifying the rational ball described by the monodromy word corresponding to the $\mathcal W$ family. Unlabeled $2-$handles carry framing $-1$.}
\end{figure}

 Our claim is that the rational ball described by the monodromy word under consideration is diffeomorphic to this one. To see this, consider the natural Kirby picture for that Lefschetz fibration, shown in Figure \ref{AI}(a) (where we have taken $p=q=r=0$ for simplicity).
Perform the indicated 1-handle slides, and then slide one 2-handle of each pair over the other as in (b). Cancelling the 1-handles that were slid with the 2-handles that were not gives Figure \ref{AI}(c), where the reader will have no trouble generalizing to arbitrary $p,q,r$. 
This last figure is identical with that obtained by ``pulling tight'' the large circles in Figure \ref{AGc}, changing $0$'s for dots, and reversing orientation.

\subsection{Graphs in family $\N$} 
\begin{figure}[t]
\includegraphics[width=5in]{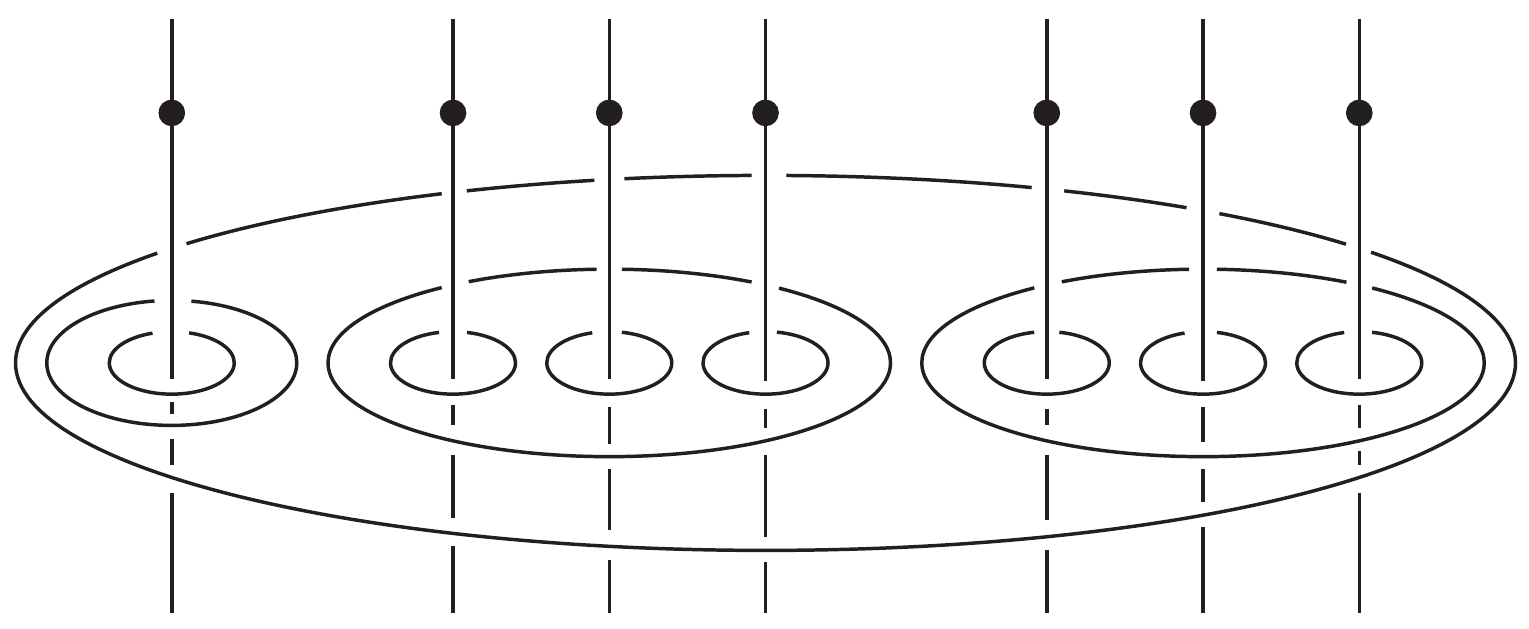}
\put(-10,40){$0$}
\caption{\label{BA}Diagram for Lefschetz fibration corresponding to base case of graphs in family $\mathcal N$. Unlabeled circles carry framing $-1$.}
\end{figure}

The construction here follows the same line as for the previous case; in particular we leave the verification that the Lefschetz fibration determined by the left side of our relation for this family is diffeomorphic to the plumbing of spheres given by the graph $\Delta_{p,q,r}$ to the reader. Note that if one draws the fiber surface as a sphere with holes rather than a disk with holes, the corresponding Kirby picture can be obtained just as before, with a dotted vertical circle for each hole, provided one encloses all of these circles by a 0-framed 2-handle (c.f. Figure \ref{BA}).

\begin{figure}[t]
\includegraphics[width=3in]{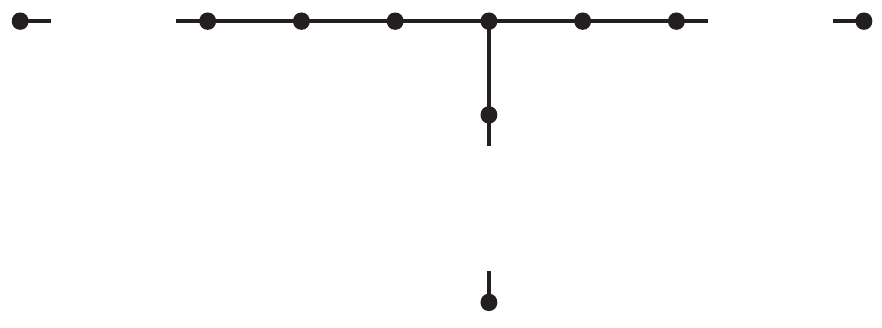}
\put(-195,73){$\cdots$}
\put(-216,71){$\underbrace{\hspace{.75in}}_{r+1}$}
\put(-162,65){${}_{-(p+2)}$}
\put(-132,65){${}_{-(q+2)}$}
\put(-99,83){${}_0$}
\put(-90,65){${}_{-(r+2)}$}
\put(-54,71){$\underbrace{\hspace{.75in}}_{q+2}$}
\put(-33,73){$\cdots$}
\put(-103,27){$\left.\begin{array}{c}\vspace*{.6in} \\ \end{array}\right\}$}
\put(-98,25){$\vdots$}
\put(-81,24){${}^{p+1}$}
\caption{\label{BB}Plumbing graph ``dual'' to $\Delta_{p,q,r}$. Unlabeled vertices carry weight $-2$.}
\end{figure}

\begin{figure}
\includegraphics[width=5in]{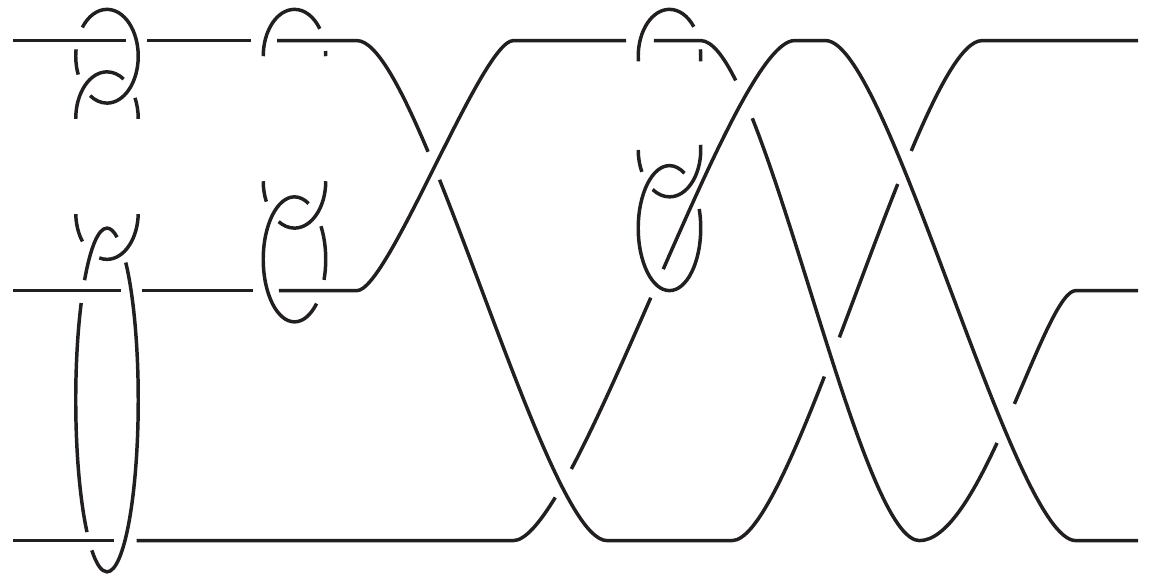}
\put(-330,125){$\vdots$}
\put(-322,130){${}_{p+1}$}
\put(-340,183){$-1$}
\put(-350,55){$0$}
\put(-272,140){$\vdots$}
\put(-264,145){${}_{r+1}$}
\put(-280,70){$-1$}
\put(-157,144){$\vdots$}
\put(-177,149){${}_{q+2}$}
\put(-153,83){$-1$}
\put(-45,175){$-(p+2)$}
\put(-45,97){$-(r+2)$}
\put(-45,3){$-(q+2)$}
\put(-75,-20){\begin{tabular}{r} $\cup$\hspace{.5ex}  3 $3$-handles\\  4-handle\end{tabular}}
\caption{\label{BC}The plumbing $\Delta'_{p,q,r}$ in a blowup of $\cee P^2$. Unlabeled circles carry framing $-2$.}
\end{figure}

\begin{figure}[t]
\includegraphics[width=4in]{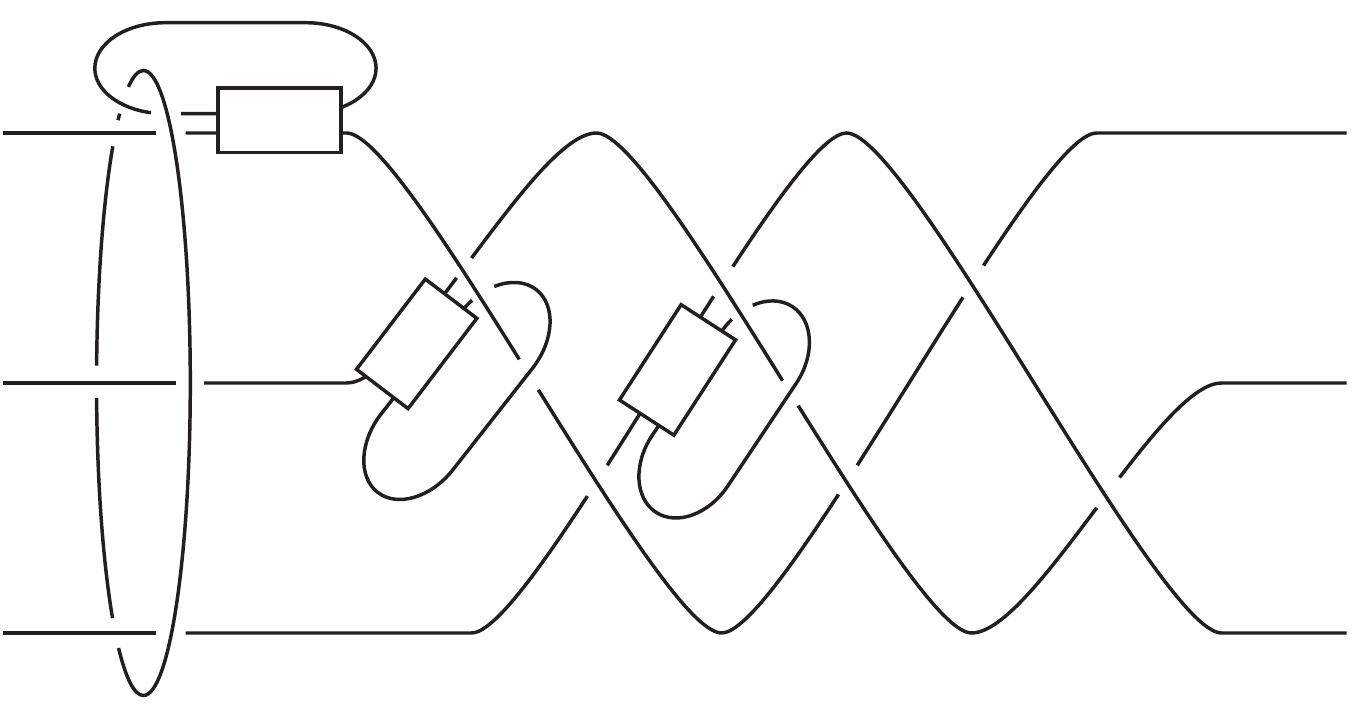}
\put(-250,150){$p+2$}
\put(-237,124){${}_{p+2}$}
\put(-205,71.5){\begin{rotate}{53}${}_{r+2}$\end{rotate}}
\put(-180,93){$r+2$}
\put(-150,66){\begin{rotate}{55}${}_{q+3}$\end{rotate}}
\put(-123,91){$q+3$}
\put(-245,90){$1$}
\put(-20,127){$1$}
\put(-20,73){$1$}
\put(-20,20){$1$}
\caption{\label{BD}After adding $0$-framed meridians to the $-1$-circles in Figure \ref{BC} and blowing down.}
\end{figure}

\begin{figure}[b]
\includegraphics[width=5in]{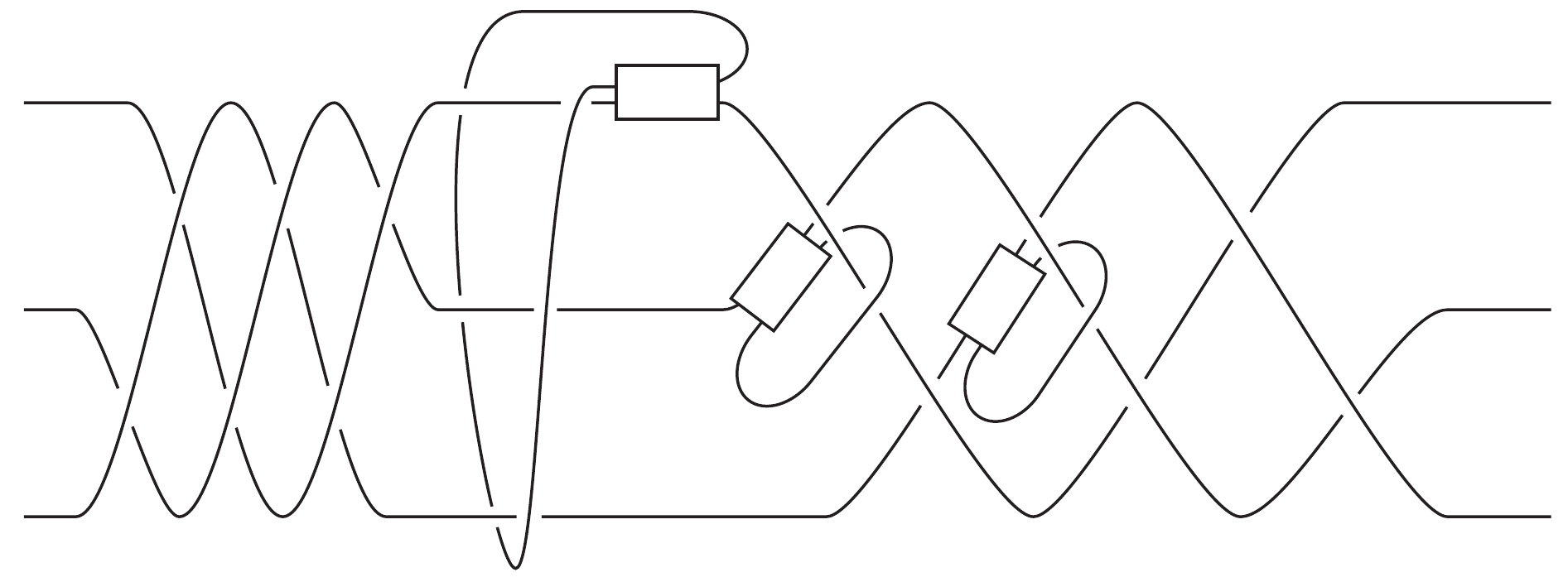}
\put(-230,137){$p+1$}
\put(-215,113){${}_{p+2}$}
\put(-186,64){\begin{rotate}{53}${}_{r+2}$\end{rotate}}
\put(-157,84){$r+2$}
\put(-136,59){\begin{rotate}{55}${}_{q+3}$\end{rotate}}
\put(-110,84){$q+3$}
\put(-20,115){$0$}
\put(-20,67){$0$}
\put(-20,20){$0$}
\caption{\label{BEa}One additional blowdown of Figure \ref{BD}. }
\end{figure}

\begin{figure}[h]
\includegraphics[width=3.8in]{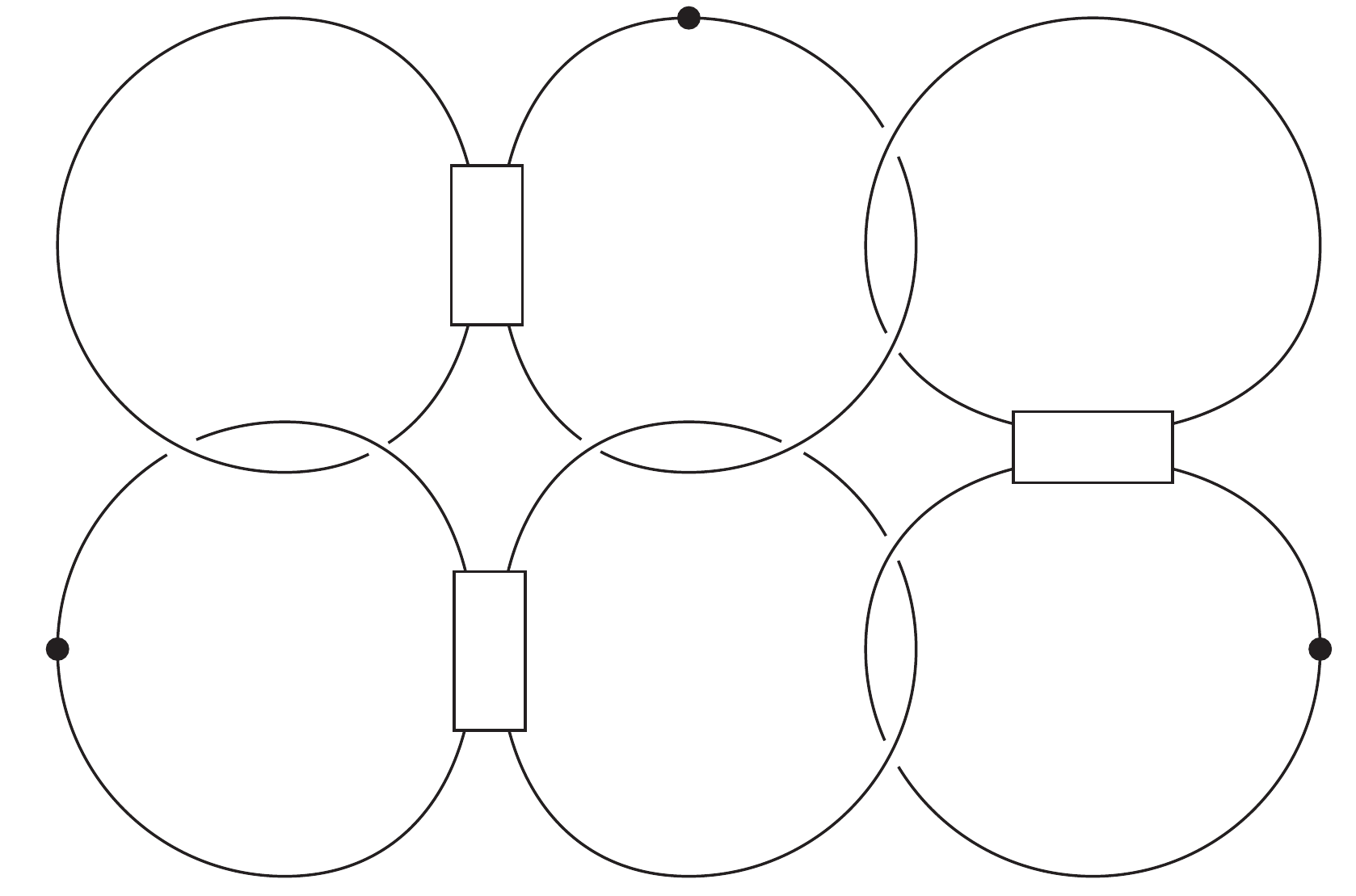}
\put(-179,119){\begin{sideways}${}_{-(r+2)}$\end{sideways}}
\put(-179,36){\begin{sideways}${}_{-(p+1)}$\end{sideways}}
\put(-240,185){$-(r+2)$}
\put(-155,-7){$-(p+1)$}
\put(-67,91){${}_{-(q+3)}$}
\put(-80,185){$-(q+3)$}
\caption{\label{BEb}Rational ball bounded by the plumbed 3-manifold described by the graph $\Delta_{p,q,r}$.}
\end{figure}

To obtain a picture of the corresponding rational homology ball, we again follow Stipsicz-Szab\'o-Wahl \cite{SSW} and begin with the dual graph $\Delta_{p,q,r}'$ shown in Figure \ref{BB}. 
It was observed in \cite{SSW} that this graph can be found spanning the rational homology in a blowup of $\cee P^2$, by blowing up a configuration of four generic lines. Beginning with Figure \ref{AE} again, this time we blow up to obtain Figure \ref{BC} which, ignoring the $-1$ circles and the $3-$ and $4-$ handles, describes the plumbing $\Delta_{p,q,r}'$.
 To find the complement of the plumbing, we add $0$-framed meridians to the $-1$ circles as before, and simplify the resulting diagram to the standard one for a connected sum of three copies of $S^1\times S^2$ (with additional framed circles corresponding to the added meridians). Adding the $0$-framed meridians and blowing back down to our original four lines gives Figure \ref{BD}. 
Blowing down the ``vertical'' 1-framed circle gives Figure \ref{BEa}, and an isotopy brings us to Figure \ref{BEb} after reversing the orientation as before.

\begin{figure}[h]
\includegraphics[width=2.8in]{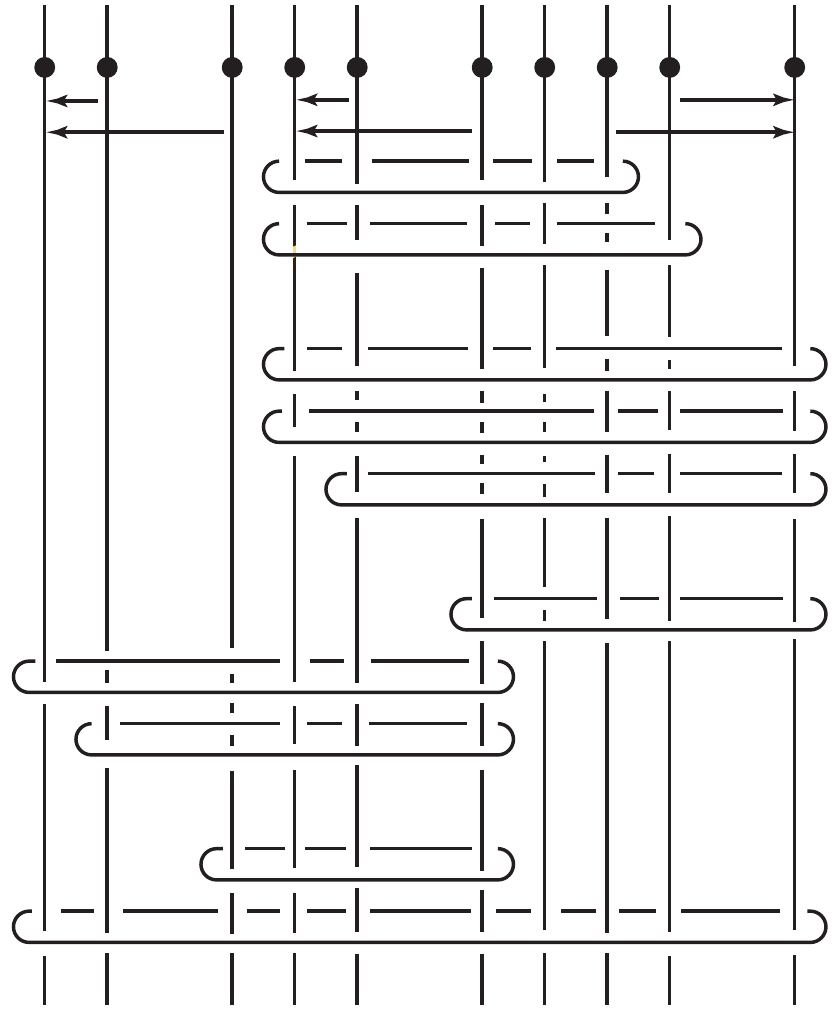}
\put(-195,4){$\underbrace{\hspace*{.72in}}_{p+1}$}
\put(-135,4){$\underbrace{\hspace*{.72in}}_{r+2}$}
\put(-59,4){$\underbrace{\hspace*{.72in}}_{q+3}$}
\put(-168,228){$\cdots$}
\put(-107.5,228){$\cdots$}
\put(-31.5,228){$\cdots$}
\put(-103,170){$\vdots$}
\put(-33,170){$\ddots$}
\put(-107,110){$\ddots$}
\put(-29,110){$\vdots$}
\put(-168,49){$\ddots$}
\put(-103,49){$\vdots$}
\put(0,20){$0$}
\caption{\label{BF}Lefschetz fibration with boundary the plumbed 3-manifold described by $\Delta_{p,q,r}$. Unlabeled circles carry framing $-1$.}
\end{figure}

\begin{figure}[h]
\includegraphics[width=2.5in]{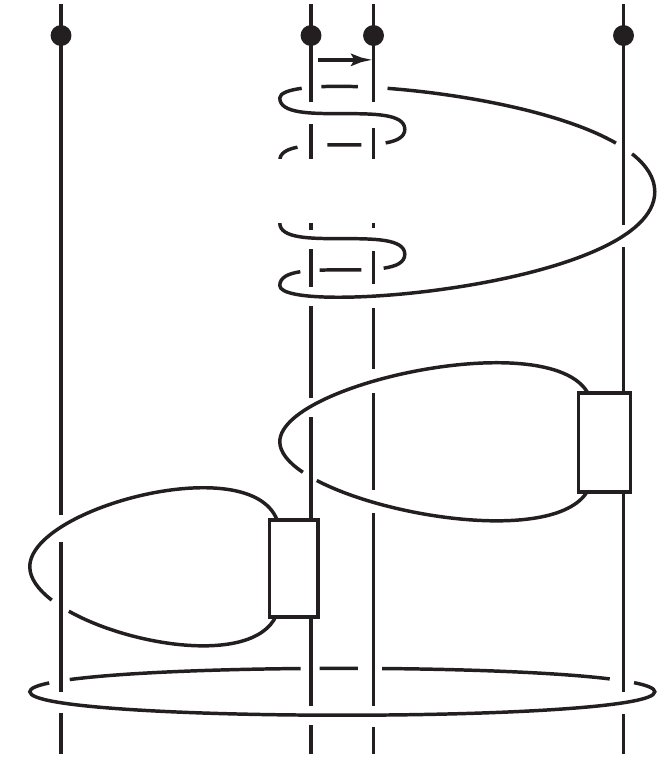}\qquad
\includegraphics[width=2.25in]{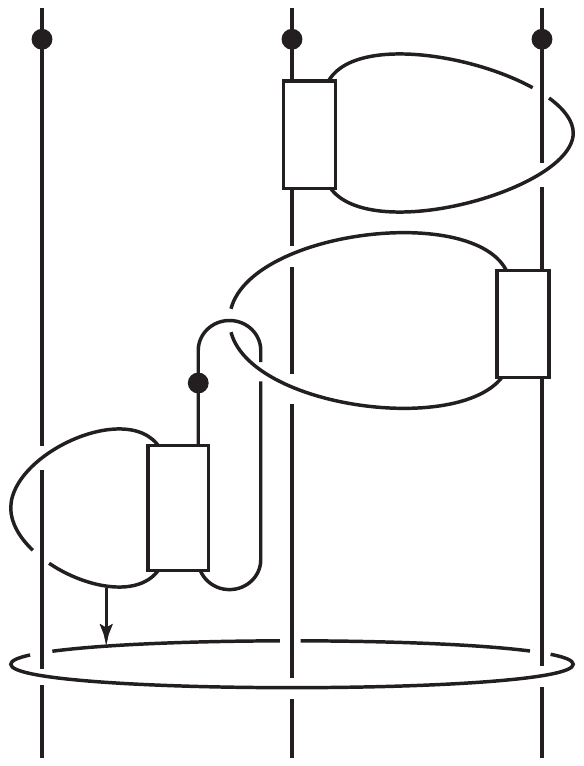}
\put(-269,154){$\vdots$}
\put(-260,158){${}_{q+3}$}
\put(-240,190){$-(q+3)$}
\put(-245,58){$-(r+2)$}
\put(-332,83){$-(p+1)$}
\put(-284,42){\begin{sideways}${}_{-(p+1)}$\end{sideways}}
\put(-200,76.5){\begin{sideways}${}_{-(r+2)}$\end{sideways}}
\put(-235,32){$0$}
\put(-145,102){${}_{-(p+1)}$}
\put(-115,60){\begin{sideways}${}_{-(p+1)}$\end{sideways}}
\put(-70,89){$-(r+2)$}
\put(-65,203){$-(q+3)$}
\put(-50,40){$0$}
\put(-22,110){\begin{sideways}${}_{-(r+2)}$\end{sideways}}
\put(-80,161){\begin{sideways}${}_{-(q+3)}$\end{sideways}}
\\\vspace{2ex}

\includegraphics[width=1.63in]{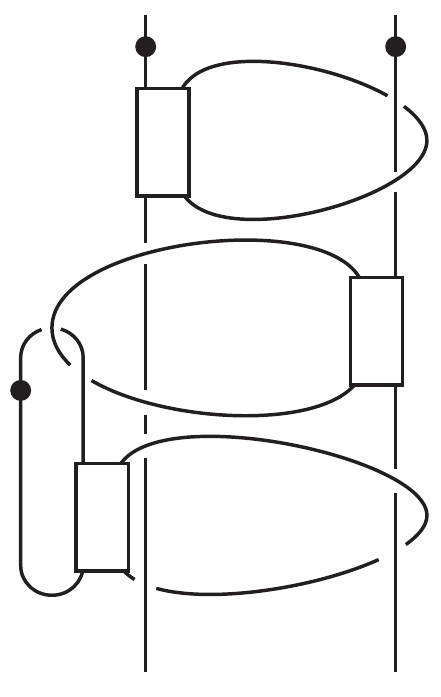}
\put(-155,200){(a)}
\put(35,200){(b)}
\put(-70,-10){(c)}
\put(-70,15){$-(p+1)$}
\put(-70,175){$-(q+3)$}
\put(-130,120){$-(r+2)$}
\put(-78,135){\begin{sideways}${}_{-(q+3)}$\end{sideways}}
\put(-94,35){\begin{sideways}${}_{-(p+1)}$\end{sideways}}
\put(-21,84){\begin{sideways}${}_{-(r+2)}$\end{sideways}}
\caption{\label{BG}Simplifying the diagram of Figure \ref{BF}.}
\end{figure}

Again, we can check that Figure \ref{BEb} describes the same rational ball as that given by the Lefschetz fibration corresponding to our monodromy word: the handle picture for the latter is shown in Figure \ref{BF}. 
Performing 1-handle slides and cancellations gives Figure \ref{BG}(c), which again is the same as the result of pulling tight the 1-handles in Figure \ref{BEa} after reversing the orientation. 
\begin{figure}[t]
\includegraphics[width=2in]{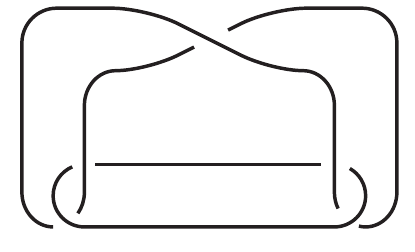}
\put(0,50){$-4$}
\put(-80,-5){$-1$}
\put(10,10){$\cup$ 3-handle}
\put(10,-2){$\cup$ 4-handle}
\caption{\label{A}A diagram for $\cptwobar$.}
\end{figure}

\newpage
\phantom{hello}
\newpage
\phantom{hello}
\newpage

\subsection{Linear plumbings} 
\begin{figure}[h]
\includegraphics[width=2in]{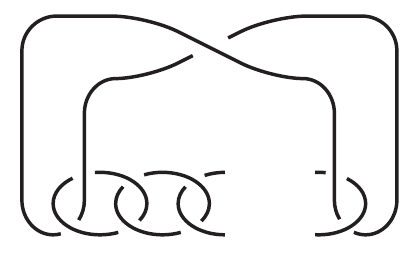}
\put(-5,60){$-(p+2)$}
\put(-57,17){$\cdots$}
\put(-120,-2){$-1$}
\put(-100,-2){$-2$}
\put(-78,-2){$-2$}
\put(-43,-2){$-2$}
\put(-123,18){$*$}
\put(10,0){$\cup$ 3-handle}
\put(10,-12){$\cup$ 4-handle}
\caption{\label{B}After $p-2$ blowups. There are $p-2$ circles with framing $-2$.}
\end{figure}

As a warmup we consider the case of $C_{p,1}$. Figure \ref{A} is a Kirby diagram for $\cptwobar$, as is easily checked. 
For $p\geq 2$, blow up one of the clasps between the $-1$ and $-4$ curves $p-2$ times to obtain Figure \ref{B}, in which we see the configuration $C_{p,1}$ embedded in $\#^{p-1} \cptwobar$. 
The complement of that configuration is given by the $-1$ framed 2-handle, together with the 3- and 4-handles, and is clearly a rational ball. Rather than pausing to derive the (well-known) Kirby diagram for this ball, we continue to the general case.

\begin{figure}[b]
\includegraphics[width=2.5in]{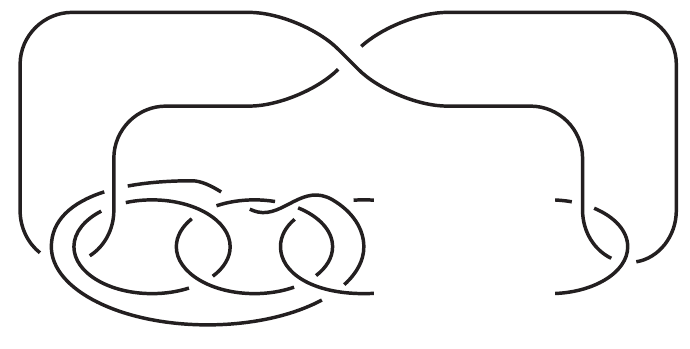}
\put(-50,95){$-(p+2)$}
\put(-150,19){${}_{-1}$}
\put(-120,19){${}_{-2}$}
\put(-170,5){$-1$}
\put(-130,24){$*$}
\put(-95,3){$-2$}
\put(-40,3){$-2$}
\put(-66,23){\Large $\cdots$}
\put(10,20){$\cup$ 2 3-handles}
\put(10,8){$\cup$ 4-handle}
\caption{\label{C}}
\end{figure}

It is well known that to construct the general $B_{p,q}$ one continues to blow up Figure \ref{B} at one or other clasp of the $-1$ curve to see a copy of $C_{p,q}$ embedded in a connected sum $\#^N \cptwobar$ and spanning the rational homology. In principle, one could carry out the procedure of the previous subsections to obtain a diagram for $B_{p,q}$ from this picture, but complications arise when the blowups alternate between clasps: in this case the curve $m$ begins to link the curve to blow down many times, and the diagram quickly becomes very complicated. To alleviate this difficulty somewhat we introduce the following device. The blowups performed to obtain Figure \ref{B} have been performed at the clasp marked with a star; we suppose that the next stage in the construction of $C_{p,q}$ calls for a blowup at the other clasp of the $-1$ curve. Introduce a cancelling pair of $2$- and $3$-handles and slide the 2-handle over the two circles between which the next blowup is to be performed (i.e., the $-1$ circle and the adjacent $-2$ circle), as indicated in Figure \ref{C}. 
Then continue the blowups using the clasp of the $-1$ curve marked with a star in that figure. Continue this procedure: each time the recipe for $C_{p,q}$ calls for the sequence of blowups to switch from one clasp to the other of the $-1$ curve, we introduce a new $2/3$-handle pair and ``save a copy'' of the $-1$ curve before proceeding. 

\begin{figure}[h]
\includegraphics[width=3in]{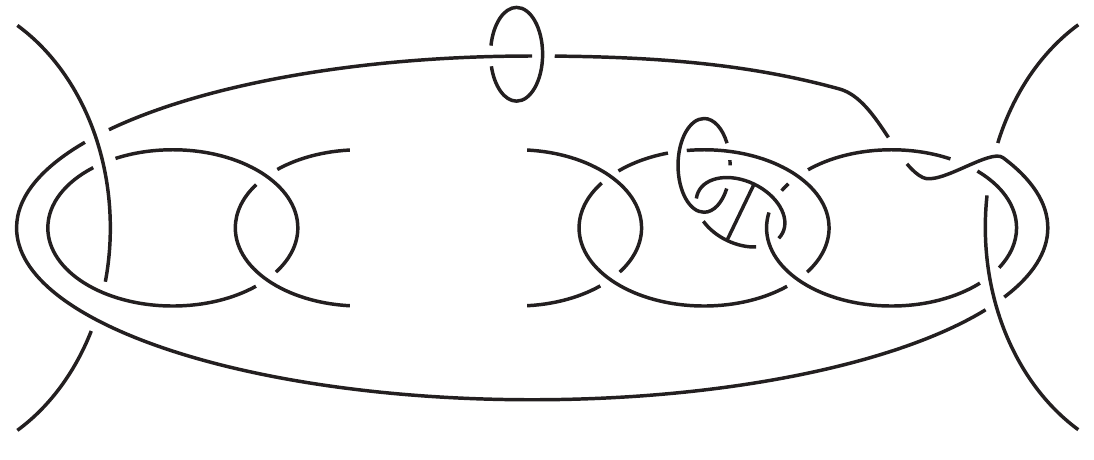}
\put(-208,82){$-c'$}
\put(-107,85){$0$}
\put(-24,82){$-2$}
\put(-87,68){$0$}
\put(-60,24){$-c$}
\put(-85,20){$-1$}
\put(-117,20){$-2$}
\put(-155,20){$-2$}
\put(-177,23){$-2$}
\put(-136,42){\Large $\cdots$}
\put(-117,0){$-1$}
\caption{\label{D}After blowing up to find the configuration $C_{p,q}$. There are $c-2$ curves with framing $-2$ in the central region.}
\end{figure}

The result of this plan is a Kirby picture for a connected sum of copies of $\cptwobar$ in which $C_{p,q}$ is visible, with some number of additional 2- and 3-handles. (The number of additional 2/3 handle pairs is one fewer than the number of continued fraction coefficients of $p^2/(pq-1)$ different from $2$.) As before, we introduce 0-framed meridians on each ``extra'' 2-handle and blow back down. After all blowdowns are complete we will be left with a 0-framed unlink together with the (now knotted and nontrivially framed) former meridians; changing the 0-framings to dots gives the Kirby diagram for $B_{p,q}$.
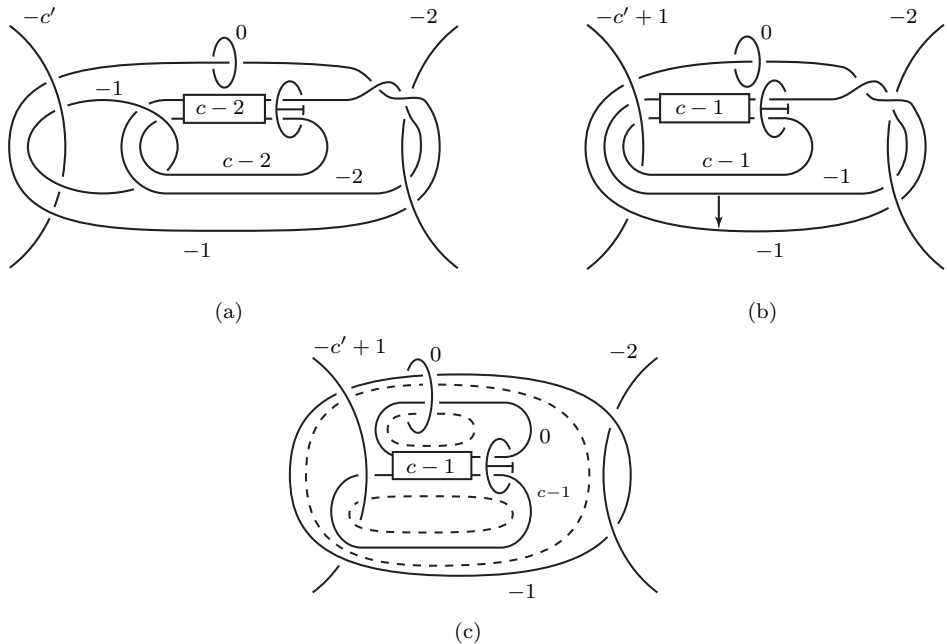
\begin{figure}[b]
\includegraphics[width=2.5in]{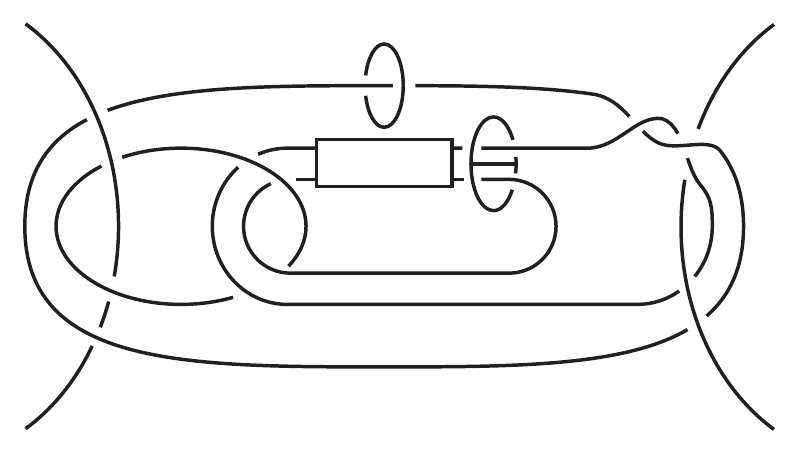}\hspace{.45in}
\includegraphics[width=2in]{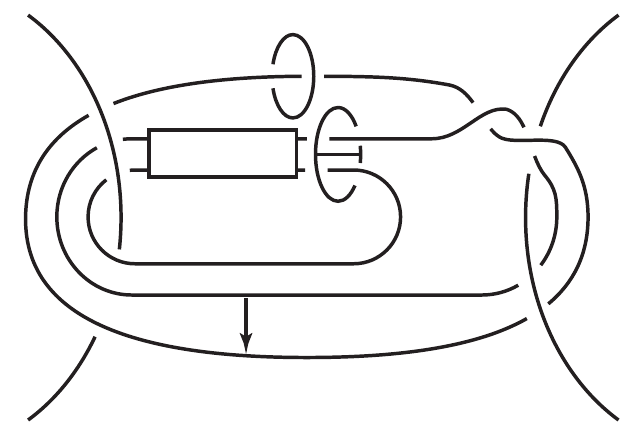}
\put(-350,98){$-c'$}
\put(-205,98){$-2$}
\put(-290,10){$-1$}
\put(-270,92){$0$}
\put(-323,71){$-1$}
\put(-285,64){$c-2$}
\put(-275,44){$c-2$}
\put(-233,38){$-2$}
\put(-135,98){$-c'+1$}
\put(-25,98){$-2$}
\put(-73,92){$0$}
\put(-105,64){$c-1$}
\put(-95,44){$c-1$}
\put(-50,38){$-1$}
\put(-75,10){$-1$}
\\\vspace{2ex}

\includegraphics[width=2in]{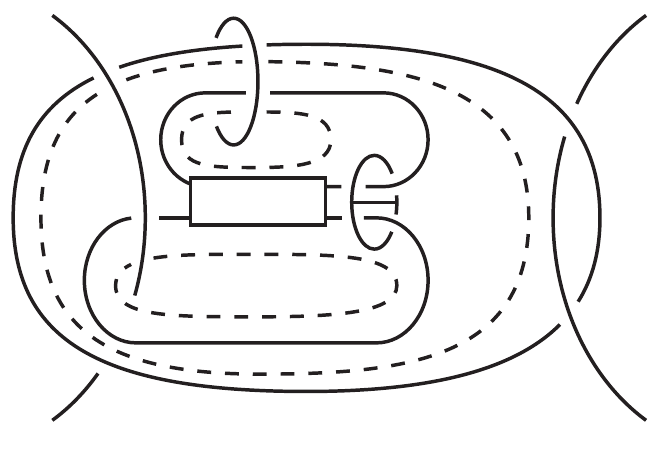}
\put(-170,110){(a)}
\put(30,110){(b)}
\put(-133,97){$-c'+1$}
\put(-89,94){$0$}
\put(-22,95){$-2$}
\put(-48,63){$0$}
\put(-49,40){${}^{c-1}$}
\put(-98,52){$c-1$}
\put(-60,5){$-1$}
\put(-80,-10){(c)}
\caption{\label{E} Simplifying Figure \ref{D}. }
\end{figure}
 To determine this diagram suppose the final sequence of blowups is as described in Figure \ref{D}; here we have drawn only the part of the diagram impinging on the last stage of the construction, which we suppose to have involved $(c-2)$ blowups for some $c\geq 3$. (In framing the outermost curves by $-c'$ and $-2$, we are also implicitly assuming that $q > 1$, i.e., this is not a diagram for $C_{p,1}$.) 
In Figure \ref{D} we have introduced $0$-framed meridians dual to the $-1$ circles, and also a $\theta$-curve that will be useful for us momentarily. Perform $(c-2)$ blowdowns to reach Figure \ref{E}(a). One additional blowdown yields Figure \ref{E}(b); sliding one $-1$ circle over the other as indicated, an isotopy then reaches Figure \ref{E}(c).

\begin{figure}[h]
\includegraphics[width=5in]{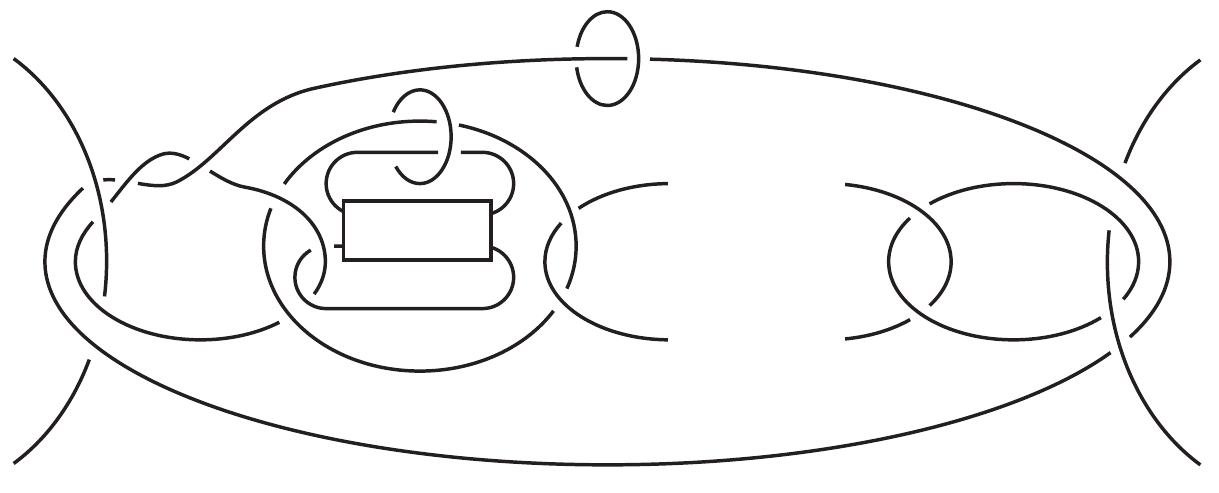}
\put(-343,123){$-2$}
\put(-32,123){$-c''$}
\put(-307,35){$-c'+1$}
\put(-245,26){$-1$}
\put(-245,45){$c-1$}
\put(-248,75){$c-1$}
\put(-207,81){$0$}
\put(-225,114){$0$}
\put(-180,34){$-2$}
\put(-147,62){\LARGE $\cdots$}
\put(-110,34){$-2$}
\put(-80,34){$-2$}
\put(-170,140){$0$}
\put(-100,125){$-1$}

\includegraphics[width=2.5in]{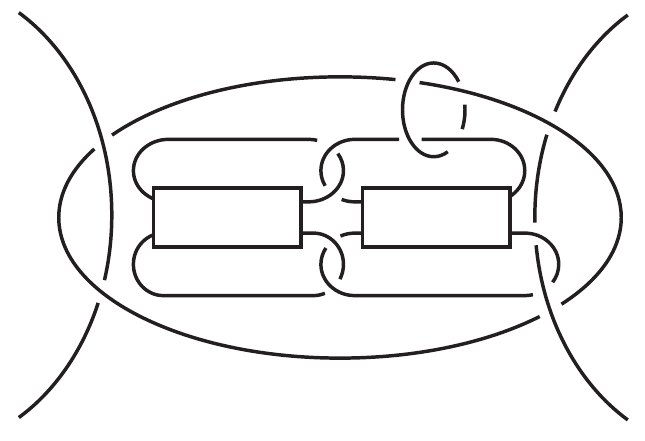}
\put(-100,130){(a)}
\put(-160,110){$-2$}
\put(-50,115){$-c''+1$}
\put(-70,108){$0$}
\put(-120,14){$-1$}
\put(-128,60){$c-1$}
\put(-120,89){$c-1$}
\put(-109,32){$0$}
\put(-80,32){$c'-2$}
\put(-71,60){$c'-2$}
\put(-80,89){$0$}
\put(-100,-5){(b)}
\caption{\label{F}(a) is Figure \ref{E} after enlarging the diagram to include the preceding sequence of blowups; there are $c'-3$ curves with framing $-2$ in the central region. Note this is the same as Figure \ref{D} after a rotation and replacing the $\theta$-curve in that diagram by the $c-1$-twisted $\theta$-link shown here. (b) shows the result of performing $c'-2$ blowdowns, following the procedure of Figure \ref{E}.}
\end{figure} 

\begin{figure}
\includegraphics[width=4in]{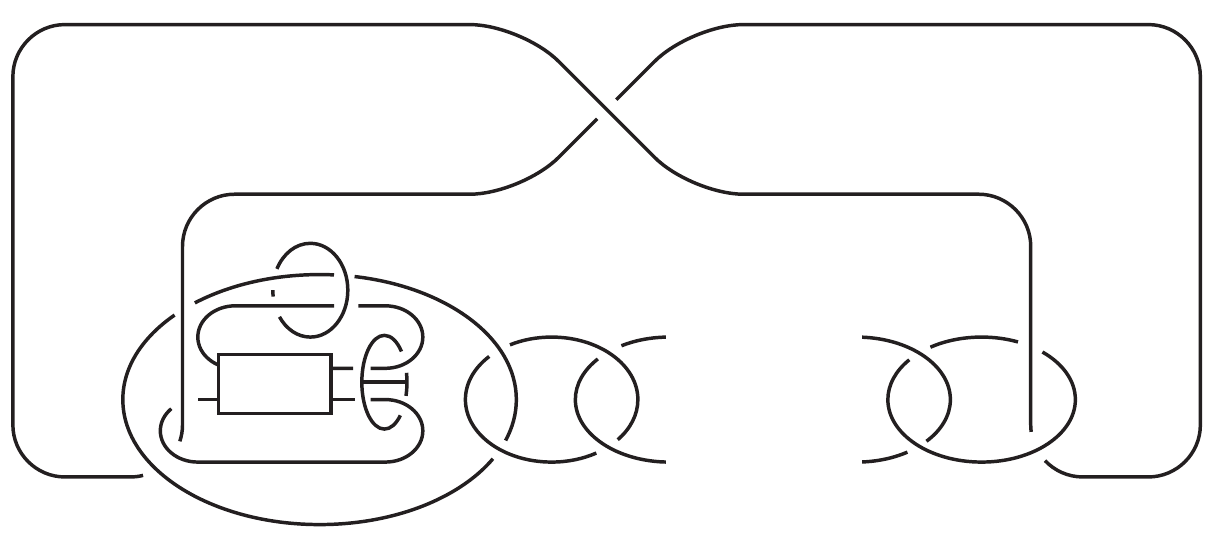}
\put(-240,125){$-(p+1)$}
\put(-235,35){ ${}_{c'-2}$}
\put(-208,68){$0$}
\put(-187,40){$0$}
\put(-220,8){$c'-2$}
\put(-210,-8){$-1$}
\put(-170,50){$-2$}
\put(-145,50){$-2$}
\put(-90,50){$-2$}
\put(-65,50){$-2$}
\put(-115,28){\Large $\cdots$}

(a)\vspace{3ex}

\includegraphics[width=2in]{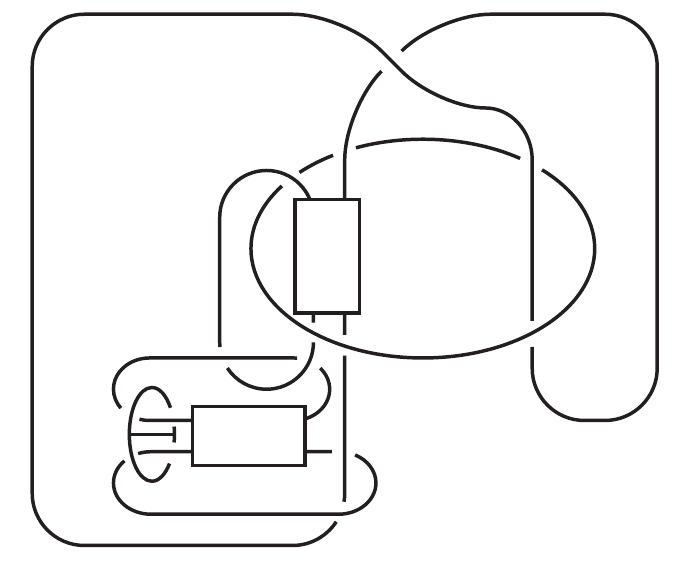}\hspace*{1in}
\includegraphics[width=1.5in]{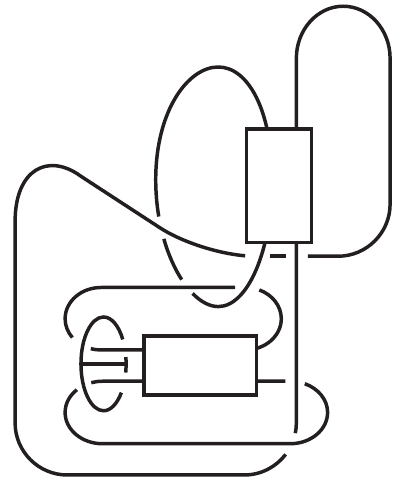}
\put(-300,125){$-4$}
\put(-284,29){${}_{c'-2}$}
\put(-262,58){\begin{sideways}\small$p-3$\end{sideways}}
\put(-300,90){$p-3$}
\put(-300,50){$0$}
\put(-245,18){${c'-2}$}
\put(-245,35){$-1$}
\put(-100,90){$0$}
\put(-70,117){$p-2$}
\put(-38,71){\begin{sideways}$p-2$\end{sideways}}
\put(-67,31){$c'-2$}
\put(-90,58){$0$}
\put(-18,20){$c'-2$}

(b)\hspace*{2.5in}(c)
\vspace{5ex}

\includegraphics[width=2.5in]{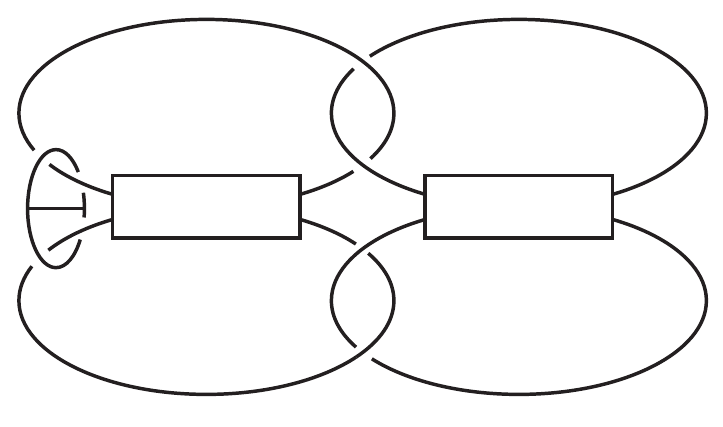}
\put(-165,100){$0$}
\put(-185,5){$c'-2$}
\put(-142,52){$c'-2$}
\put(-60,52){$p-1$}
\put(-20,100){$p-2$}
\put(-20,5){$0$}

(d)
\caption{\label{G}The endgame of simplification. Diagram (a) is obtained from Figure \ref{C} by replacing the configuration of three curves enclosing the $*$ in that figure by the diagram of Figure \ref{E}(a), and blowing down once following the moves of Figure \ref{E}(b) and (c). The framings follow the conventions of Figure \ref{C} (for the large curve, which we recall has been involved in one blowdown) and Figure \ref{F} (in that we assume at least one more $\theta$-link appears where indicated). In (a) there are  $p-3$ curves with framing $-2$; (b) results from (a) by $p-3$ blowdowns. An additional blowdown gives (c), which is isotopic to (d).}
\end{figure}

The iterative step is now reasonably clear. Think of the linking $0$- and $c-1$-framed circles as lying in a genus 2 handlebody, which is a neighborhood of a $\theta$-curve linking the $-c'+1$ curve and the $0$-framed meridian (the boundary of this neighborhood is indicated by the dashed lines in Figure \ref{E}(c)). To continue our blowdowns with the $-1$-circle followed by the $-2$-circles to the right in Figure \ref{E}(c), we can return with a $180^\circ$ rotation to Figure \ref{D}, with the $\theta$-curve in that figure replaced by the ``$\theta$-link'' in Figure \ref{E}(c). Note, however, that one blowdown in this next stage has already been performed, corresponding to the transition from Figure \ref{E}(a) to (b). The result of this iteration is shown in Figure \ref{F}, and the process clearly continues with another ``$\theta$-link'' added in the center of Figure \ref{F}(b) each time the sequence of blowdowns switches from one clasp of the $-1$ curve to the other. (Note that while the framing and number of twists appearing in the first $\theta$-link differs from the framing appearing in the original plumbing by one, for successive $\theta$-links this number is adjusted by 2: in other words the next $\theta$-link to appear in this process will have framing and twisting number both equal to $c''-2$.) 

The final stage of the construction is shown in Figure \ref{G}, which depicts the situation after blowing down Figure \ref{C} once and applying the construction of Figure \ref{E}, with the additional $\theta$-links indicated.

To state the end result we need only determine the framings and numbers of twists each $\theta$-link receives.

\begin{figure}[t]
\includegraphics[width=5in]{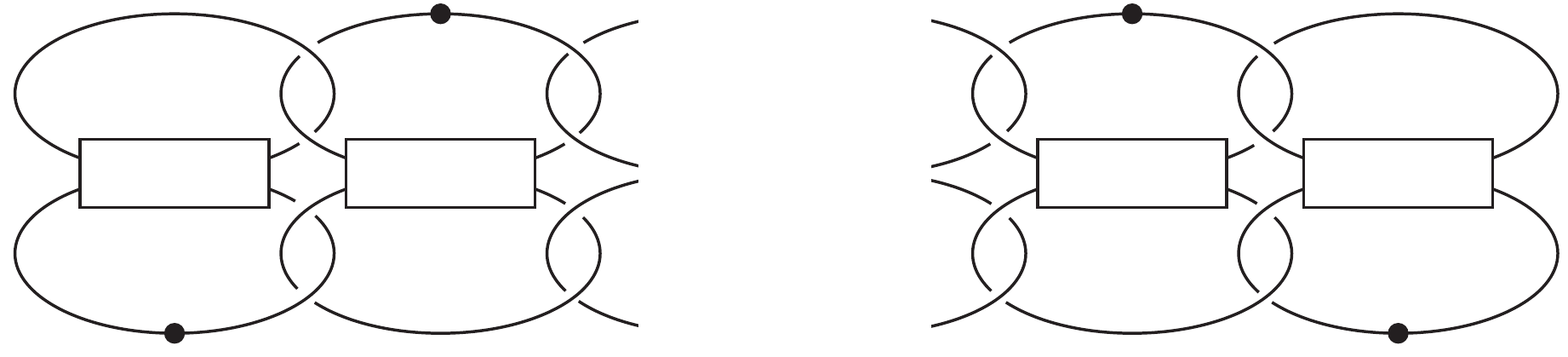}
\put(-340,83){$-(c_1-4)$}
\put(-337,42){${}_{-(c_1-3)}$}
\put(-276,42){${}_{-(c_2-2)}$}
\put(-279,-3){$-(c_2-2)$}
\put(-120,42){${}_{-(c_{\ell-1}-2)}$}
\put(-123,-3){$-(c_{\ell-1}-2)$}
\put(-55,42){${}_{-(c_\ell -1)}$}
\put(-58,83){$-(c_\ell -1)$}
\put(-188,38){\LARGE $\cdots$}
\caption{\label{H} The rational ball bounded by the plumbed 3-manifold corresponding to the graph $\Gamma_{p,q}$. The numbers $c_1,\ldots,c_\ell$ are the absolute values of the integers different from $-2$ appearing in the continued fraction expansion $-p^2/(pq-1) = [-b_1,\ldots,-b_n]$, written in a particular order. Namely, if we construct the continued fraction by applications of the moves (a) and (b) in Theorem \ref{linplumbthm}, we list the $c_j$ in the order in which they appear in the construction. Equivalently, if the expansion is obtained by applying (a) $x_1$ times, then (b) $x_2$ times, and so on, then $x_1 = c_1-4$, while $x_j =  c_j-2$ for $j = 2,\ldots,\ell$.}
\end{figure}

\begin{theorem}\label{linplumbthm}Let $[b_1,\ldots,b_k]$ be the continued fraction expansion of $p^2/(pq-1)$ (where each $b_i\geq 2$), so that $[-b_1,\ldots,-b_k]$ can be obtained from the length-one expansion $[-4]$ by repeated application of the operations
\begin{enumerate}
\item[a)] $[-b_1,\ldots,-b_k] \mapsto [-b_1-1,-b_2,\ldots, -b_k, -2]$
\item[b)] $[-b_1,\ldots,-b_k] \mapsto [-2, -b_1,\ldots, -b_{k-1}, -b_k-1]$.
\end{enumerate}
Let $\{c_1,\ldots, c_\ell\}$ be the sequence of integers different from $2$ appearing among the $\{b_1,\ldots, b_k\}$, written in the order in which they arise during the application of (a) and (b) above. Then a Kirby diagram for a rational homology ball $B_{p,q}$ with boundary $L(p^2,pq-1)$ is given by Figure \ref{H}.
\end{theorem}

Indeed, the numbers $c_1-4, c_2-2,\ldots, c_\ell -2$ are the number of times each operation (a) or (b) is applied consecutively in the construction of the plumbing. From the description of $B_{p,q}$ before, we see that both the number of twists and the framing on one side of the typical $\theta$-link is just this number $c_i-2$; the exceptions are at the beginning, where the $\theta$-link gets one extra blowdown (contributing both a twist and a framing adjustment to the $\ell$-th link), and at the end, where the last blowdown (indicated in Figure \ref{G}) results in one extra twist but no adjustment to the framing. (Note that in the degenerate case that only one number $c$ obtains, both effects arise: in this case the diagram is a single $\theta$ link with $-(c-2)$ twists and framing $-(c-3)$.) To obtain Figure \ref{H}, we have performed all the blowdowns, changed the 0-framed unknots to dotted circles, and reversed the orientation to account for the fact that our construction involved inverting a manifold.

Turning to the monodromy relation, the fact that one side of the relation obtained previously describes a linear plumbing of spheres is an easy exercise.
The verification that the rational ball described above is diffeomorphic to the one given by the monodromy word follows much the same lines as the previous arguments, and we leave it to the interested reader. 

\section{An example of substitution}

Here we construct a family of rational blowdowns along linear chains $C_{p,1}$ via monodromy substitution using the daisy relation.

\begin{figure}[t]
\includegraphics{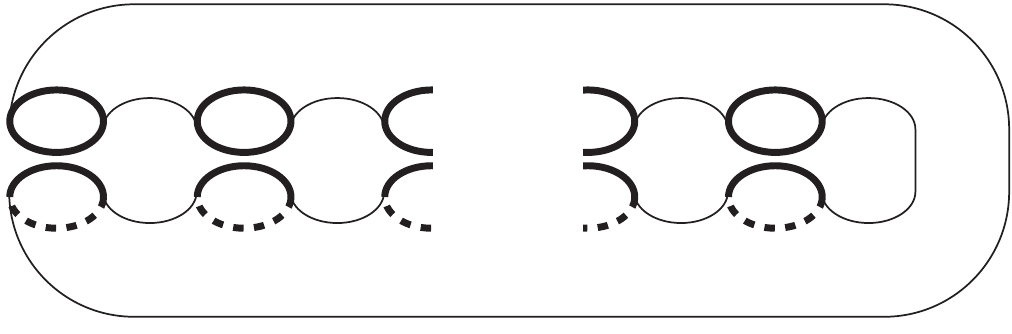} 
\put(-275,75){$d_1$}
\put(-230,75){$d_2$}
\put(-70,75){$d_{k+1}$}
\put(-275,20){$d$}
\put(-230,20){$d_{p+1}$}
\put(-70,20){$d_{k+2}$}
\put(-153,60){\Large $\cdots$}
\put(-153,35){\Large $\cdots$}
\put(0,20){$S_{p+2}$}
\vspace{2ex}\\
\includegraphics{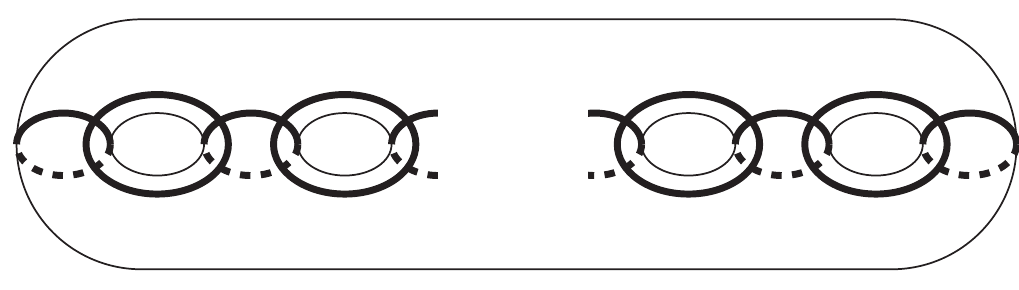}
\put(-278,55){$c_1$}
\put(-253,20){$c_2$}
\put(-228,55){$c_3$}
\put(-200,20){$c_4$}
\put(-153,40){\Large $\cdots$}
\put(-105,20){$c_{2g-2}$}
\put(-80,55){$c_{2g-1}$}
\put(-48,20){$c_{2g}$}
\put(-28,55){$c_{2g+1}$}
\put(0,20){$\Sigma_g$}
\caption{\label{exfig}}
\end{figure}

Suppose that $p$ is even and $p=2k\, (k\geq 1)$. 
We put $g=k+1$ and take a chain $c_1,c_2,\ldots ,c_{2g+1}$ of length $2g+1$ 
on $\Sigma_g$ as in Figure \ref{exfig}. 
We choose an embedding $f:S_{p+2}\rightarrow\Sigma_g$ of a planar surface with $p+2$ boundary components such that 
$f(d_i)=f(d_{p-i+3})=c_{2i-1}\, (i=1,2,\ldots ,k+1)$, where we set $d_{p+2}:=d$ 
(see Figure \ref{exfig}). 
The embedding $f$ maps the left-hand side $d_1d_2\cdots d_pd_{p+1}d^{p-1}$ 
of the daisy relation to the word 
\[
c_1^{2k}c_3^2c_5^2\cdots c_{2k+1}^2=c_1^{2g-2}c_3^2c_5^2\cdots c_{2g-1}^2.
\]
Our goal is to find a word in positive Dehn twists representing the identity in the mapping class group $\Mod(\Sigma_g)$, containing the above as a subword.

We first consider the word
\[
C_{\rm I}:=(c_1c_2c_3\cdots c_{2g}c_{2g+1})^{2g+2},
\]
which represents the identity in $\Mod(\Sigma_g)$.
Applying braid relations, we have the following sequence of equivalences. 
\begin{eqnarray*}
C_{\rm I}&=&(c_1c_2c_3\cdots c_{2g}c_{2g+1})^{2g+2} \\
&\equiv & c_1\cdot c_2c_1\cdot c_3c_2\cdot c_4c_3\cdot c_5c_4\cdots c_{2g}c_{2g-1}\cdot 
c_{2g+1}c_{2g}\cdot c_{2g+1}\cdot (c_1c_2c_3\cdots c_{2g}c_{2g+1})^{2g} \\
&\equiv & c_1\cdot c_2c_1\cdot c_3c_2c_1\cdot c_4c_3c_2\cdot c_5c_4c_3\cdot
c_6c_5c_4\cdots c_{2g+1}c_{2g}c_{2g-1}\cdot c_{2g+1}c_{2g}c_{2g+1} \\
&& \cdot (c_1c_2c_3\cdots c_{2g}c_{2g+1})^{2g-1} \\
&\equiv & c_1c_3\cdot c_2c_3\cdot c_1c_2\cdot c_4c_5\cdot c_3c_4\cdot c_2c_3\cdot 
c_6c_5c_4\cdots c_{2g+1}c_{2g}c_{2g-1}\cdot c_{2g+1}c_{2g}c_{2g+1} \\
&& \cdot c_1c_2c_1\cdot c_3c_4\cdots c_{2g}c_{2g+1}\cdot c_2c_3\cdots c_{2g}c_{2g+1}
\cdot (c_1c_2c_3\cdots c_{2g}c_{2g+1})^{2g-3} \\
&\equiv & c_1c_3c_5\cdot c_2c_3\cdot c_1c_2\cdot c_4c_5\cdot c_3c_4\cdot c_2c_3\cdot 
c_6c_5c_4\cdots c_{2g+1}c_{2g}c_{2g-1}\cdot c_{2g+1}c_{2g}c_{2g+1} \\
&& \cdot c_1c_2c_3c_4\cdots c_{2g}c_{2g+1}\cdot c_2c_3\cdots c_{2g}c_{2g+1}
\cdot (c_1c_2c_3\cdots c_{2g}c_{2g+1})^{2g-3} \\
&\equiv & c_1c_3c_5\cdot c_2c_3\cdot c_1c_2\cdot c_4c_5\cdot c_3c_4\cdot c_2c_3\cdot 
c_6c_7\cdot c_5c_6\cdot c_4c_5 \\
&& \cdot c_8c_7c_6\cdots c_{2g+1}c_{2g}c_{2g-1}\cdot c_{2g+1}c_{2g}c_{2g+1} 
\cdot c_1c_2c_3c_4\cdots c_{2g}c_{2g+1}\cdot c_2c_3\cdots c_{2g}c_{2g+1} \\
&& \cdot c_1c_2c_1\cdot c_3c_4\cdots c_{2g}c_{2g+1}\cdot c_2c_3\cdots c_{2g}c_{2g+1}
\cdot (c_1c_2c_3\cdots c_{2g}c_{2g+1})^{2g-5} \\
&\equiv & c_1c_3c_5c_7\cdot c_2c_3\cdot c_1c_2\cdot c_4c_5\cdot c_3c_4\cdot c_2c_3\cdot 
c_6c_7\cdot c_5c_6\cdot c_4c_5 \\
&& \cdot c_8c_7c_6\cdots c_{2g+1}c_{2g}c_{2g-1}\cdot c_{2g+1}c_{2g}c_{2g+1} \\
&& \cdot (c_1c_2c_3c_4\cdots c_{2g}c_{2g+1}\cdot c_2c_3\cdots c_{2g}c_{2g+1})^2 
\cdot (c_1c_2c_3\cdots c_{2g}c_{2g+1})^{2g-5} \\
&\equiv & \cdots \\
&\equiv & c_1c_3c_5\cdots c_{2g-1}\cdot 
c_2c_3\cdot c_1c_2\cdot c_4c_5\cdot c_3c_4\cdot c_2c_3\cdot c_6c_7\cdot c_5c_6\cdot c_4c_5 \\
&& \cdots c_{2g-2}c_{2g-1}\cdot c_{2g-3}c_{2g-2}\cdot c_{2g-4}c_{2g-3}\cdot 
c_{2g}c_{2g-1}c_{2g-2}\cdot c_{2g+1}c_{2g}c_{2g-1}\cdot c_{2g+1}c_{2g}c_{2g+1} \\
&& \cdot (c_1c_2c_3c_4\cdots c_{2g}c_{2g+1}\cdot c_2c_3\cdots c_{2g}c_{2g+1})^{g-2} 
\cdot (c_1c_2c_3\cdots c_{2g}c_{2g+1})^3 \\
&= & c_1c_3c_5\cdots c_{2g-1}\cdot W, 
\end{eqnarray*}
where we put 
\begin{eqnarray*}
W&:= & 
c_2c_3\cdot c_1c_2\cdot c_4c_3c_2\cdot c_5c_4c_3\cdot c_6c_5c_4\cdot c_7c_6c_5 \\
&& \cdots c_{2g-2}c_{2g-3}c_{2g-4}\cdot c_{2g-1}c_{2g-2}c_{2g-3}\cdot 
c_{2g}c_{2g-1}c_{2g-2}\cdot c_{2g+1}c_{2g}c_{2g-1}\cdot c_{2g+1}c_{2g}c_{2g+1} \\
&& \cdot (c_1c_2c_3c_4\cdots c_{2g}c_{2g+1}\cdot c_2c_3\cdots c_{2g}c_{2g+1})^{g-2} 
\cdot (c_1c_2c_3\cdots c_{2g}c_{2g+1})^3. \\
\end{eqnarray*}
We then make the positive relator 
\[
\bar{C}_{\rm I}:=(c_{2g+1}c_{2g}\cdots c_3c_2c_1)^{2g+2} 
\]
into another relator $\bar{W}\cdot c_{2g-1}\cdots c_5c_3c_1$ in a similar way, 
where we put 
\begin{eqnarray*}
\bar{W}&:= & 
(c_{2g+1}c_{2g}\cdots c_3c_2c_1)^3\cdot 
(c_{2g+1}c_{2g}\cdots c_3c_2\cdot c_{2g+1}c_{2g}\cdots c_4c_3c_2c_1)^{g-2} \\
&& \cdot c_{2g+1}c_{2g}c_{2g+1}\cdot c_{2g-1}c_{2g}c_{2g+1}\cdot c_{2g-2}c_{2g-1}c_{2g} 
\cdot c_{2g-3}c_{2g-2}c_{2g-1}\cdot c_{2g-4}c_{2g-3}c_{2g-2} \\
&& \cdots c_5c_6c_7\cdot c_4c_5c_6\cdot c_3c_4c_5\cdot c_2c_3c_4\cdot c_2c_1 
\cdot c_3c_2. 
\end{eqnarray*}
The hyperelliptic relator 
\[
I^2:=(c_1c_2c_3\cdots c_{2g}c_{2g+1}^2c_{2g}\cdots c_3c_2c_1)^2 
\]
is conjugate to the relator $c_1^2\cdot U$, where we put 
\[
U:=c_2c_3\cdots c_{2g}c_{2g+1}^2c_{2g}\cdots c_3c_2c_1^2c_2c_3\cdots 
c_{2g}c_{2g+1}^2c_{2g}\cdots c_3c_2. 
\]
Gathering these relators, we obtain a factorization of the identity containing the left-hand side of the daisy relation: 
\begin{eqnarray*}
\bar{W}\cdot c_{2g-1}\cdots c_5c_3c_1\cdot c_1c_3c_5\cdots c_{2g-1}\cdot W
\cdot (c_1^2\cdot U)^{g-2}\hspace*{-1.5in} \\
&\equiv & \bar{W}\cdot c_{2g-1}\cdots c_5c_3c_1\cdot c_1^2\cdot 
c_1c_3c_5\cdots c_{2g-1}\cdot W\cdot (c_1^2\cdot U)^{g-3}\cdot U \\
&\equiv & \bar{W}\cdot c_{2g-1}\cdots c_5c_3c_1\cdot c_1^4\cdot 
c_1c_3c_5\cdots c_{2g-1}\cdot W\cdot (c_1^2\cdot U)^{g-4}\cdot U^2 \\
&\equiv & \cdots \\
&\equiv & \bar{W}\cdot c_{2g-1}\cdots c_5c_3c_1\cdot c_1^{2g-4}\cdot 
c_1c_3c_5\cdots c_{2g-1}\cdot W\cdot U^{g-2} \\
&\equiv & \bar{W}\cdot c_1^{2g-2}c_3^2c_5^2\cdots c_{2g-1}^2\cdot W\cdot U^{g-2} \\
&=: & \varrho.
\end{eqnarray*}

The word $\varrho$ is a positive relator and consists of non-separating 
simple closed curves. The corresponding Lefschetz fibration $X(\varrho)$ is a fiber sum of 
two copies of $X(C_{\rm I})$ and $g-2$ copies of $X(I^2)$. 
The signature of $X(\varrho)$ is $-4(2g-1)(g+1)$ (cf. \cite{endo}) 
and the Euler characteristic 
of $X(\varrho)$ is $4g(4g-1)$. This manifold is simply-connected and minimal 
(see \cite{usher}). We claim that $X(\varrho)$ is not spin, from which it follows that $X(\varrho)$ is homeomorphic but not diffeomorphic to 
$\#(4g^2-4g+1)\Bbb{CP}^2\#(12g^2-3)\overline{\Bbb{CP}}^2$. To see that $X(\varrho)$ is not spin, recall that it is isomorphic to a fiber sum of 
$X(C_{\rm I}^2)$ and $X(I^{2(g-2)})$ by construction. 
Since $X(C_{\rm I}^2)$ is isomorphic to $X(I^{2(g+1)})$ 
(\cite{auroux} Lemma 3.4, cf. \cite{endo1} Lemma 4.1), 
we can show that $X(\varrho)$ is isomorphic to $X(I^{2(2g-1)})$, 
a fiber sum of $2g-1$ copies of $X(I^2)$. 
Since $X(I^2)$ has a section with square $-1$, 
$X(I^{2(2g-1)})$ has a section with square $-(2g-1)$. 
Hence the intersection form of $X(I^{2(2g-1)})$ is odd, 
which implies $X(I^{2(2g-1)})$ is not spin.


If we apply a daisy substitution to $\varrho$ 
we obtain a new positive relator $\varrho'$, and 
the corresponding Lefschetz fibration $X(\varrho')$ 
is a rational blowdown of $X(\varrho)$ along a configuration $C_{p,1}$
The signature of $X(\varrho')$ is $-8g^2-2g+1$ 
and the Euler characteristic 
of $X(\varrho')$ is $16g^2-6g+3$. 
This manifold is simply-connected, non-spin, and symplectic. 
Hence $X(\varrho')$ is homeomorphic but not diffeomorphic to 
$\#(4g^2-4g+1)\cee P^2\#(12g^2-2g)\cptwobar$ (cf. \cite{GS}, sec. 2.4). 

Likewise, if we apply a daisy substitution to the product $\varrho I^4$ of $\varrho$ with $I^4$ 
we obtain a positive relator $\varrho' I^4$, and
the corresponding Lefschetz fibration $X(\varrho' I^4)$ 
is a rational blowdown of $X(\varrho I^4)$ along a $C_{p,1}$. 
The three manifolds $X(\varrho I^4)$, $X(\varrho' I^4)\# (2g-3)\cptwobar$, 
and $\#(4g^2+1)\cee P^2\#(12g^2+12g+5)\cptwobar$ are 
homeomorphic but mutually non-diffeomorphic (cf. \cite{EG}, Theorem 4.1). 

{\bf Remark.} 
We can construct similar $\varrho$ and $\varrho'$ for odd $p>1$.

\affiliationone{Hisaaki Endo\\ Department of Mathematics\\
Graduate School of Science\\
Osaka University\\
Toyonaka, Osaka 560-0043\\
JAPAN\\
\email{endo@math.sci.osaka-u.ac.jp}}
\affiliationtwo{Thomas E. Mark\\ Department of Mathematics\\ University of Virginia\\ PO Box 400137\\ Charlottesville, VA 22904\\ USA \\ \email{tmark@virginia.edu}}
\affiliationthree{Jeremy Van Horn-Morris \\American Institute of Mathematics\\
360 Portage Ave\\
Palo Alto, CA 94306-2244\\
USA\\
\email{jvanhorn@aimath.org}}

\end{document}